\documentclass[graybox]{svmult}

\usepackage{type1cm}        
\usepackage{makeidx}        
\usepackage{graphicx}       
\usepackage{multicol}        
\usepackage[bottom]{footmisc}

\usepackage{newtxtext} 

\makeindex            
\usepackage{amsmath}
\usepackage{amsfonts}
\usepackage{microtype}
\usepackage{cleveref}
\usepackage{authblk}

\crefname{figure}{\textbf{Figure}}{\textbf{Figures}}
\crefname{section}{\textbf{Section}}{\textbf{Sections}}
\crefname{equation}{\textbf{Equation}}{\textbf{Equations}}
\crefname{remark}{\textbf{Remark}}{\textbf{Remarks}}
\crefname{theorem}{\textbf{Theorem}}{\textbf{Theorems}}
\crefname{definition}{\textbf{Definition}}{\textbf{Definitions}}
\crefname{lemma}{\textbf{Lemma}}{\textbf{Lemmas}}
\crefname{proposition}{\textbf{Proposition}}{\textbf{Propositions}}
\crefname{table}{\textbf{Table}}{\textbf{Tables}}

\graphicspath{{./figs/}}

\usepackage{tikz}
\usepackage{subcaption}
\usepackage{tabularx}
\usepackage{enumitem}

\usepackage{xcolor}
\usepackage{multirow}

\DeclareMathOperator{\proj}{proj}

\usepackage{appendix}
\usepackage{amsmath}
\usepackage{amsthm}
\usepackage{amsfonts} 
\usepackage{tabularx}
\usepackage{tikz-cd}

\newcommand {\mm}[1] {\ifmmode{#1}\else{\mbox{\(#1\)}}\fi}

\newcommand{\Nspace}        {\mm{{\mathbb N}}}
\newcommand{\Rspace}        {\mm{{\mathbb R}}}

\newcommand{\Fcal}        {\mm{\mathcal F}}

\newcommand{\Hcal}        {\mm{\mathcal H}}
\newcommand{\Kcal}        {\mm{\mathcal K}}
\newcommand{\Lcal}        {\mm{\mathcal L}}
\newcommand{\Ncal}        {\mm{\mathcal N}}

\newcommand{\Qcal}        {\mm{\mathcal Q}}
\newcommand{\Rcal}        {\mm{\mathcal R}}
\newcommand{\Scal}        {\mm{\mathcal S}}

\newcommand{\Hgroup}        {\mm{{\mathsf H}}}
\newcommand{\Cgroup}        {\mm{{\mathsf C}}}
\newcommand{\rank}{\mathrm{rank}}

\newcommand{\OR}{\kappa}
\newcommand{\Ric}{\operatorname{Ric}}

\newcommand{\dime}[1]       {\mm{\rm dim\,}{#1}}

\newcommand{\Cech}        {\mm{{\mathrm{\check{C}ech}}}}
\newcommand{\RVR}        {\mm{{\mathrm{\hat{\Rcal}}}}}
\newcommand{\VR}        {\mm{{\mathrm{VR}}}}
\newcommand{\Vor}        {\mm{{\mathrm{Vor}}}}
\newcommand{\Del}        {\mm{{\mathrm{Del}}}}
\newcommand{\Wrap}        {\mm{{\mathrm{Wrap}}}}
\newcommand{\DelCech}        {\mm{{\mathrm{Del\check{C}ech}}}}

\newcommand{\diam}        {\mm{{\mathrm{diam}}}}
\newcommand{\image}        {\mm{\mathrm{im}}}

\newcommand{\para}[1]{\vspace{2mm}\noindent{\textbf{#1}}}

\newcommand\TODO[3]{\hbox to 0pt{\textcolor{#1}{$^\bullet$}}\marginpar{\footnotesize \textcolor{#1}{\begin{flushleft}#2: #3\end{flushleft}}}}

\title*{Finding the Cores of Higher Graphs Using Geometric and Topological Means: A Survey}
\titlerunning{Finding the Cores of Higher Graphs}
\author{Inés García-Redondo\orcidID{0000-0001-8340-4235}, 
Claudia Landi\orcidID{0000-0001-8725-4844}, 
Sarah Percival\orcidID{0000-0003-1024-4618},
Anda Skeja\orcidID{0000-0002-5707-5420},
Bei Wang\orcidID{0000-0002-9240-0700},
and Ling Zhou\orcidID{0000-0001-6655-5162}}
\authorrunning{García-Redondo et al.}
\institute{Inés García-Redondo \at London School of Geometry and Number Theory, Imperial College London, UK.\\
\email{i.garcia-redondo22@imperial.ac.uk}
\and Claudia Landi \at DISMI, University of Modena and Reggio Emilia, Italy.\\
\email{clandi@unimore.it}
\and Sarah Percival \at Department of Mathematics and Statistics, University of New Mexico, USA.\\ \email{spercival@unm.edu}
\and Anda Skeja \at Department of Mathematics, Uppsala University, Sweden.\\
\email{anda.skeja@math.uu.se}
\and Bei Wang \at School of Computing, Scientific Computing and Imaging Institute, University of Utah, USA.\\
\email{beiwang@sci.utah.edu}
\and Ling Zhou \at Department of Mathematics, Duke University, USA.\\
\email{ling.zhou@duke.edu}
}

\begin{document}

\maketitle 

\abstract{
In this survey, we explore recent literature on finding the cores of higher graphs using geometric and topological means. We study graphs, hypergraphs, and simplicial complexes, all of which are models of higher graphs. We study the notion of a core, which is a minimalist representation of a higher graph that retains its geometric or topological information. We focus on geometric and topological methods based on discrete curvatures, effective resistance, and persistent homology. We aim to connect tools from graph theory, discrete geometry, and computational topology to inspire new research on the simplification of higher graphs. \newline

}

\section{Introduction}
\label{sec:introduction}
Real-world networks are often modeled as graphs,  where edges capture pairwise interactions between entities represented by vertices.
However, in many social and biological scenarios, pairwise interactions are mere reductions from actual multiway interactions, which can be captured more adequately by hypergraphs and simplicial complexes. 
In studying pathogenic viral responses, multiple genes (vertices) coordinate with each other to form biological pathways (hyperedges)~\cite{FengHeathJefferson2021}. 
While investigating consensus formation, opinions of individuals (vertices) in a society are not formed in isolation but rather in large groups (hyperedges)~\cite{GolovinMolterKuehn2024}. 
We study graphs, hypergraphs, and simplicial complexes, all of which are considered as forms of higher graphs\footnote{Spivak~\cite{Spivak2009} considered these objects to belong to the categories of higher graphs and studied them in a unifying framework. We simply borrow the notion of higher graphs in this paper to encompass the above objects of interest.}.  
On the other hand, although a simplicial complex can be viewed as a special case of a hypergraph—specifically, as a downward-closed hypergraph—these two structures possess distinct algebraic properties (e.g., homology theories) and are therefore treated separately in this context.

In this paper, we survey the recent literature on finding the \emph{cores} of higher graphs using geometric and topological means.
Informally, the core of a higher graph is a minimalist representation that retains its geometric or topological information. 
Here, the concept of a core is more general than the \emph{$k$-core} of a graph~\cite{Seidman1983}, which is a maximal connected subgraph in which each  vertex has a degree of at least $k$.
The concept of a core also goes beyond the \emph{$(k,q)$-core} of a hypergraph~\cite{AhmedBatageljFu2007}, which is a maximal subgraph in which each vertex has a hypergraph degree of at least $k$ and each hyperedge contains at least $q$ vertices~\cite{LeeGohLee2023}. 
The notion of a core has been used in~\cite{MemoliOkutan2021} to refer to a minimalistic simplicial filtration that retains the persistent homology information.
Several notions from network science are closely related to this concept. 
\emph{Coarse geometry} of networks preserves geometric or topological properties while ignoring small-scale features~\cite{WeberJostSaucan2018}.  
A \emph{network backbone} retains short-range interactions while disregarding long-range ones~\cite{BensonGleichLeskovec2016}. 
The study of cores is relevant to graph reduction methods, including sparsification, coarsening, and condensation, some of which focus on preserving geometric or topological information within the reduced representations; see~\cite{HashemiGongNi2024,LiuSafaviDighe2018} for surveys. 

We focus on geometric and topological methods in extracting the cores of higher graphs, based on discrete curvatures, effective resistance, and persistent homology.
We do not survey these three classes of concepts in general, but rather focus on using them for core extraction. 
We aim to connect tools from graph theory, discrete geometry, and computational topology to inspire new research on the simplification of higher graphs. 

\para{Paper classification.} 
The annotation of each paper is guided primarily by three types of higher graphs, namely, \textbf{graphs}, \textbf{hypergraphs}, and \textbf{simplicial complexes}. 
It is guided secondarily by certain geometric and topological methods used to preserve the cores of higher graphs, which, admittedly, partially reflect our own biases in constructing this survey: \textbf{discrete curvatures}, which include Forman-Ricci, Ollivier-Ricci, and resistance curvatures; \textbf{effective resistance}, which arises from various notions of Laplacians; and \textbf{homology}, which includes simplicial homology and persistent homology. For hypergraphs and simplicial complexes, we also survey geometric methods from \textbf{percolation} and \textbf{spectral clustering}. 

\para{Connection to existing surveys.}
There is a vast literature on graph reduction, also known as graph summarization or graph simplification;  see~\cite{ChenSaadZhang2022,ChenYeVedula2023,HashemiGongNi2024,InterdonatoMagnaniPerna2020,LiuSafaviDighe2018, VonLuxburg2007} for surveys. 
\emph{Graph sparsification} approximates a given graph by a sparse graph with a subset of vertices and/or edges~\cite{ChenYeVedula2023}. 
\emph{Graph coarsening} groups nodes into super-nodes and aggregates the intergroup edges into super-edges~\cite{ChenSaadZhang2022}.  
\emph{Graph condensation} condenses a graph by synthesizing a smaller graph with comparable performance for graph neural networks (GNNs)~\cite{JinZhaoZhang2022};  see~\cite{GaoYuJiang2024} for a survey. 
Different from existing surveys, we focus on higher graph sparsification and coarsening techniques based on tools in combinatorics (e.g., Laplacians, discrete curvatures) and computational topology (e.g., persistent homology).

Both graph coarsening and graph clustering involve grouping nodes, but with different goals: graph coarsening seeks to reduce the size of the graph while maintaining its essential properties, whereas graph clustering aims to identify meaningful communities or patterns within the data. The literature on spectral clustering of graphs is extensive (e.g., \cite{ShiMalik2000,NgJordanWeiss2001,MeilaShi2001}), focusing on detecting communities among nodes by clustering the eigenvectors of specific Laplacians; for a comprehensive survey, see \cite{VonLuxburg2007}. 
For the purpose of this survey, we exclude spectral graph clustering due to its long history and popularity; instead, we include recent advancements on the spectral clusterings of hypergraphs and simplicial complexes as they create new research opportunities.  

In network science, percolation theory has received considerable attention, describing the behaviors of a network (e.g, giant clusters, cluster distribution) when nodes or edges are randomly designated either occupied or unoccupied~\cite{LiLiuLu2021}. We exclude network percolation (including an in-depth discussion on $k$-cores~\cite{KongShiWu2019}) from this survey;  see~\cite{Saberi2015,LeeKahngCho2018,LiLiuLu2021} for recent reviews. Instead, we pay special attention to the recent theory of percolation on hypergraphs (e.g.~\cite{BianconiDor2024}), whereas the theory of percolation on simplicial complexes remains scarce and does not yet lead to the computation of
well-defined cores (e.g.~\cite{ZhaoLiPeng2022,ZhaoLiPeng2022b}). 

Finally, there are a number of surveys on the combinatorics and topological tools we utilize: see~\cite{Merris1994} for a classic review of graph Laplacians and~\cite{WeiWei2023} for a recent survey on persistence Laplacians derived from computational topology;~\cite{EdelsbrunnerHarer2008,EdelsbrunnerHarer2022} for introductory texts on persistent homology. 
 
\section{Technical Background on Discrete Curvatures}
\label{sec:curvature}

To study the geometry of graphs, discrete curvatures are  natural measures that quantify the local geometry around nodes and edges of a graph. 
In this survey, we are interested in finding the cores of higher graphs using discrete curvatures. 
Therefore, in this section we review various notions of discrete curvatures for higher graphs, followed by  discussions of their connections to core-findings in graphs (\cref{sec:graph-cores}), hypergraphs (\cref{sec:hypergraph-cores}), and simplicial complexes (\cref{sec:sc-cores}), respectively.   

The concept of Ricci curvature, originating from differential geometry, has seen extensive research dedicated to its discretization, for application in graphs and similar structures~\cite{Stone1976,ChowLuo2003,Morgan2005,LottVillani2009,GuSaucan2013,BonciocatSturm2009,JinKimGu2007,AlsingMcDonaldMiller2011}. 
We begin with the discretization of Ricci curvature, namely, the Forman-Ricci curvature (\cref{sec:FR}), the Ollivier-Ricci curvature (\cref{sec:OR}), and their variants. 
We also discuss the more recent notion of resistance curvature (\cref{sec:resistance-curvature}), which is related to the effective resistance and the sparsification of graphs and simplicial complexes. 

\subsection{Forman-Ricci Curvature}
\label{sec:FR}

\para{History.} One of the simplest discrete curvatures to compute in graphs is the Forman-Ricci curvature, originally formulated for general CW-complexes by Forman \cite{Forman2003} employing a reinterpretation of the Bochner-Weitzenb\"{o}ck formula.  For an oriented, compact Riemannian manifold,  
the Hodge Laplacian can be expressed as a sum of the Bochner Laplacian, i.e.~$\nabla^* \nabla$, where $\nabla$ is the Levi-Civita connection of the manifold and $\nabla^*$ its adjoint,  and an endomorphism of the bundle of differential forms only  involving the Ricci curvature of the manifold \cite[Theorem 9.4.1]{petersen2016riemannian}. 

Forman adapted this approach to a combinatorial setting and demonstrated that, for a weighted CW complex, a similar decomposition could be achieved, and thus defined a Ricci curvature from the corresponding term within it \cite{Forman2003}. Sreejith et al.~\cite{SreejithMohanrajJost2016} introduced the Forman-Ricci curvature to network analysis, particularized for undirected, possibly weighted graphs (seen as 1-dimensional simplicial complexes). Forman-Ricci curvature was later extended to possibly directed hypergraphs by Leal et al.~\cite{LealRestrepoStadler2021}. All these notions are reviewed and unified by Eidi et al.~\cite{EidiFarzamLeal2020}.

\para{Definitions for graphs and simplicial complexes.} We begin by introducing the notion of curvature introduced by Forman \cite{Forman2003} applied to the case of \emph{simplicial complexes} (instead of CW complexes), as they constitute one of our primary objects of study.

\begin{definition}[\cite{Forman2003}]
\label{def:forman-ricci-sc}
    For a weighted simplicial complex \(S\) and any dimension $p$, the \emph{Forman-Ricci curvature} is a function
    \[\Fcal_p : S_p \to \Rspace\]
    defined for each \(p\)-simplex \(\alpha \in S_p\) as
    \begin{equation}
    \label{eq:forman-ricci-weighted-sc}
        \Fcal_p(\alpha) = w_\alpha \left( \left[\sum_{\beta>\alpha} \dfrac{w_\alpha}{w_\beta} + \sum_{\gamma<\alpha} \dfrac{w_\gamma}{w_\alpha} \right]\nonumber - \sum_{\Tilde{\alpha} \neq \alpha} \left| \sum_{\beta>\alpha,\, \Tilde{\alpha}} \dfrac{\sqrt{w_\alpha w_{\Tilde{\alpha}}}}{w_\beta} - \sum_{\gamma<\alpha,\, \Tilde{\alpha}} \dfrac{w_\gamma}{\sqrt{w_\alpha w_{\Tilde{\alpha}}}}\right|\right),
    \end{equation}
    where \(\alpha < \beta\) denotes that \(\alpha\) is contained in the boundary of the \((p+1)\)-simplex \(\beta\). In other words, \(\alpha\) is a face of \(\beta\) or \(\beta\) is a coface of \(\alpha\) of relative dimension 1.  
\end{definition}
This definition can be particularized for unweighted simplicial complexes as
\begin{align}
\label{eq:forman-ricci-unweighted-sc}
\Fcal_p(\alpha)  = & \# \{(p+1)\text{-cells }\beta > \alpha\} + \# \{(p-1)\text{-cells }\gamma < \alpha\} \nonumber \\
& - \# \{\text{parallel neighbours of }\alpha\}, 
\end{align}
where a parallel neighbour of a \(p\)-simplex \(\alpha\) is defined as another \(p\)-simplex \(\Tilde{\alpha}\) sharing only one face or one coface of relative dimension 1 with \(\alpha\). 

Assuming an undirected, unweighted graph \(G = (V,E)\), upon substituting \(p=1\) into \cref{def:forman-ricci-sc}, one derives the following curvature expression \cite{SreejithMohanrajJost2016, EidiFarzamLeal2020} for an edge \(e = (u,v)\), 
\begin{equation}
\label{eq:forman-ricci-networks-degree}
    \Fcal(e) = 2 - \deg (e).
\end{equation}
Here, \(\deg(e)\) denotes the number of parallel neighbours of an edge $e$, that is, the number of edges sharing a vertex with \(e\). Note that one can also write \(\deg(e) = \deg_u(e) + \deg_v(e) = \deg(u) + \deg(v) - 2\) where \(\deg(v)\) is the degree of a vertex $v$ and \(\deg_v(e)\) is the number of edges sharing the vertex \(v\) with \(e = (u,v)\). This relation on degrees gives the following equivalent expression for the Forman-Ricci curvature in terms of the degrees of the vertices: \(\Fcal(e) = 4 - \deg(u) - \deg(v)\).
From this notion, Sreejith et al.~\cite{SreejithMohanrajJost2016} also defined a \emph{node Forman curvature} as
\begin{equation}
    \Fcal(v) = \dfrac{1}{\deg(v)} \sum_{e_v } \Fcal(e_v), 
\end{equation}
where \(e_v\) denotes an edge incident on the node \(v\).  

In~\cite{SreejithJostSaucan2017, SaucanSreejithVivek-Ananth2019}, the Forman-Ricci curvature for nodes and edges was extended to directed graphs. Borrowing notations from \cite{EidiFarzamLeal2020}, let \(G=(V, E)\) be an unweighted, directed graph, where we denote the directed edge going from \(u\) (the tail) to \(v\) (the head) as \(e =  [u,v]\). For an edge \(e = [u, v]\), we call the \emph{input} of \(e\) to the set of all edges having \(u\) as their head, and analogously, the output of \(e\) is the set of all edges having \(v\) as its tail. The degree of the input \(\deg_{{in}}(e)\) is the number of elements in the input, and similarly we can define the degree of the output \(\deg_{{out}}(e)\). Using these notions, we can define the Forman-Ricci curvature of the edge as \(\Fcal(e) = 2 - \deg_{{in}}(e) - \deg_{{out}}(e) \).

\para{Definitions for hypergraphs.} 
From the formulation above, the Forman-Ricci curvature can be very naturally extended to (directed) hypergraphs \(H = (V, E)\) \cite{LealRestrepoStadler2021} where $V$ is the set of vertices and $E$ a set of ordered pairs of subsets of $V$, called hyperedges. For $e = (e_1, e_2) \in E$, let $e_1$ be the tail of the hyperedge, its nodes being input nodes, and $e_2$ the head of the hyperedge, composed by output nodes. We now define $$\deg_{{in}}(e) = \# \{\text{hyperedges with an input node of }e\text{ as their head}\}$$ and $$\deg_{{out}}(e) = \# \{\text{hyperedges with an output node of }e\text{ as their tail}\}.$$ 
Notice that another hyperedge \(\tilde{e}\) from \(H\) might have several input nodes of \(e\) as their head, so that the formula above is actually counting the hyperedges with multiplicity given by this number of ``shared'' nodes; and similarly for the output nodes of \(e\) and the hyperedges having them as their tail.  
\begin{definition}[\cite{LealRestrepoStadler2021}]
    Given a directed hypergraph \(H = (V, E)\) the Forman-Ricci curvature of an hyperedge $e = (e_1, e_2) \in E$ is 
    \[\Fcal(e) = |e_1| + |e_2| - \deg_{{in}}(e) - \deg_{{out}}(e).\]
\end{definition}

Note that \cite{LealRestrepoStadler2021} also generalized the definition of Forman–Ricci curvature to the setting of weighted undirected hypergraphs. Their formulation is motivated by the following expression for the curvature of an edge \(e = (e_1, e_2)\) in a weighted graph:
\[
\Fcal(e) 
= w_e \,\Bigg( \bigg( \frac{w_{e_1}}{w_e} - \sum_{e'\ni e_1,\, e'\neq e} \frac{w_{e_1}}{\sqrt{w_e w_{e'}}}\bigg) + \bigg( \frac{w_{e_2}}{w_e} - \sum_{e'\ni e_2,\, e'\neq e} \frac{w_{e_2}}{\sqrt{w_e w_{e'}}}\bigg) \Bigg).
\]
\begin{definition}[\cite{LealRestrepoStadler2021}]
Given a weighted undirected hypergraph \(H = (V, E)\), the Forman–Ricci curvature of a hyperedge \(e \in E\) is defined as
\[
\Fcal(e) = w_e \sum_{v \in e} \bigg( \frac{w_v}{w_e} - \sum_{e'\ni v,\, e'\neq e} \frac{w_v}{\sqrt{w_e w_{e'}}} \bigg).
\]
\end{definition}

\para{Properties.}
From its definition, one of the key advantages of the Forman-Ricci curvature becomes apparent: it has a direct, inexpensive computation for networks, in contrast to other notions of curvature that will be introduced subsequently. 
It is also clear from \cref{eq:forman-ricci-networks-degree} that the Forman-Ricci curvature of an edge should be very negative (have high absolute value) when the vertices defining the edge are well-connected within the network.
To investigate in more detail which features of a graph this function captures, Sreejith et al.~\cite{SreejithMohanrajJost2016} examined the Forman-Ricci curvature for nodes and edges in both theoretical and real-world networks, including transportation and protein interaction networks. They studied their distributions and investigated correlations between Forman-Ricci curvature and centrality/connectivity measures, such as degree, clustering coefficient, and betweenness centrality. As expected, strong positive correlations with \emph{degree} are found in random and small-world networks, whereas weaker negative correlations exist in scale-free networks. 
Their study also suggested an observable correlation between Forman-Ricci  curvature and the degree of nodes, as a function of degree assortativity. 
Negative correlations were discovered between Forman-Ricci curvature and various centrality measures in different network types, such as \emph{betweenness centrality}, which measures the number of shortest paths passing through an edge and is expensive to compute. Interestingly, the \emph{clustering coefficient}---a quantity typically used as a notion of curvature in network---showed no correlation with Forman-Ricci curvature. Finally, the paper examined the impact of removing nodes based on their Forman-Ricci curvature on network connectivity, finding that such removal accelerates network disintegration compared to random or clustering coefficient-based node removal. In a follow-up paper, Sreejith et al.~\cite{SreejithJostSaucan2017} studied similar notions for directed graphs, arriving at the same conclusions.
\subsection{Ollivier-Ricci Curvature}
\label{sec:OR}

\para{History.} Recall that in Riemannian geometry, the Ricci curvature is positive when ``small spheres are closer (in transportation distance) than their centers are''~\cite{ProkhorenkovaSamosvatHoorn2020}.  
Inspired by this, Ollivier used a probability measure $m_x$ that depends on $x$ serving as an analog for the sphere $S_x$ centered at $x$. 
Rather than corresponding points between two close spheres $S_x$ and $S_y$, the transportation distance between measures is utilized. 
Ollivier introduced a notion of \emph{coarse Ricci curvature} \cite{Ollivier2007,Ollivier2009}, 
now commonly referred to as the \emph{Ollivier-Ricci curvature}, of Markov chains valid on arbitrary metric spaces, such as graphs. 

With a modification of Ollivier's definition, Lin et al.~\cite{LinLuYau2011} introduced another notion of discrete curvature on graphs. 
By incorporating gradients of the graph Laplacian, M{\"u}nch and Wojciechowski \cite[Theorems 2.1 and 2.6]{MunchWojciechowski2019} 
obtained a limit-free formula of the definition of Lin et al. 
In these definitions, the Ollivier-Ricci curvature of an edge is a local quantity that depends on the degrees of its endpoints, which may result in a lack of robustness of the Ollivier-Ricci curvature for sparse networks. 
For this reason, instead of structural neighbourhoods, Gosztolai and Arnaudon~\cite{GosztolaiArnaudon2021} considered distributions generated by diffusion processes across scales.

\begin{figure}[!ht]
    \centering
    \includegraphics[width=0.75\linewidth]{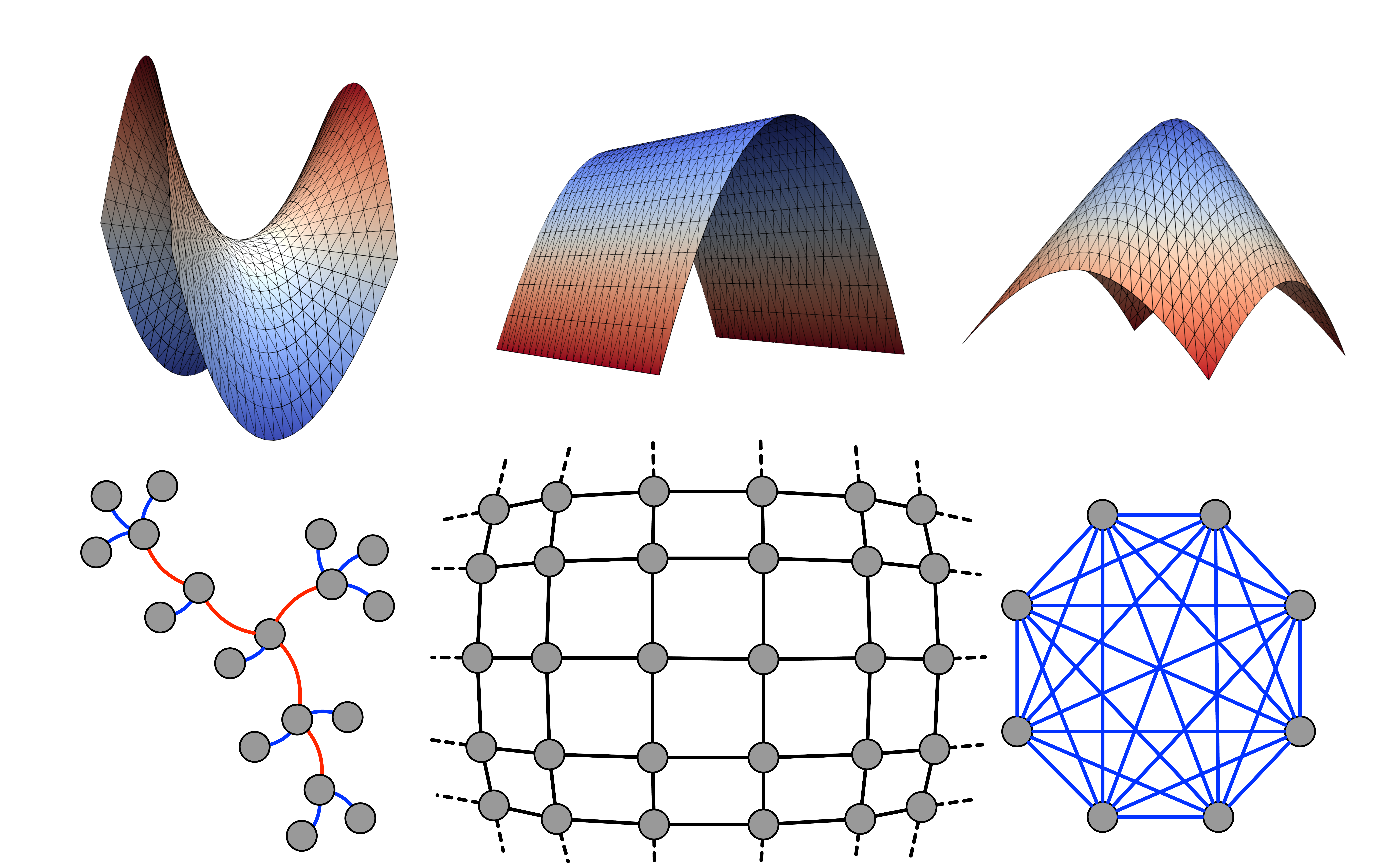}
    \caption{An illustration of Ricci curvatures for surfaces (top) and Ollivier-Ricci  curvatures (\cref{eqn:mx}) for graphs (bottom). Top: from left to right, surfaces of negative, zero, and positive Ricci curvature, respectively. Bottom: red for negative Ollivier-Ricci curvature, black for zero Ollivier-Ricci curvature, and blue for positive Ollivier-Ricci curvature. Bottom left: a tree graph with negative Ollivier-Ricci curvature everywhere except the edges connecting the leaves. Bottom middle: an infinitely sized grid graph with all edges of zero Ollivier-Ricci curvature. Bottom right: a complete graph with all edges of positive Ollivier-Ricci curvature. Notice the difference with respect to the Forman-Ricci curvature (\cref{eq:forman-ricci-networks-degree}) which would be negative in all the edges of the graphs on the bottom. Image reproduced from~\cite[Figure 4]{NiLinLuo2019} with modifications, licensed under CC BY 4.0.~\url{https://creativecommons.org/licenses/by/4.0/}.}  
    \label{fig:OR}
\end{figure}

In a further extension, Eidi and Jost \cite{EidiJost2020} generalized the Ollivier-Ricci curvature to (possibly directed and/or weighted) hypergraphs. More recently, Coupette et al.~\cite{CoupetteDalleigerRieck2023} provided a unified framework for several generalizations of Ollivier-Ricci curvatures \cite{AsoodehGaoEvans2018,EidiJost2020,LealEidiJose2020,Banerjee2021} on hypergraphs and their computations. 
Recently, Yamada \cite{Yamada2023} extended the Ollivier-Ricci curvature to simplicial complexes. 

\para{Definitions for graphs.} 
Let $(X, d)$ be a metric space and $\nu_1,\nu_2$ probability measures on $X$. 
An $L_1$ transportation distance between $\nu_1$  and $\nu_2$ is 
\[
W_1(\nu_1,\nu_2):= \inf_{\xi \in \Pi(\nu_1,\nu_2)} \int_{(x,y)\in X\times X} d(x, y) \, \mathrm{d}\xi(x, y), 
\]
where $\Pi(\nu_1,\nu_2)$ is the set of measures on $X \times X$ projecting onto $\nu_1$ and $\nu_2$, respectively.

A \emph{random} walk $m$ on $X$ is defined as a family of probability measures $m_x(\cdot)$ on $X$ for each $x\in X$, satisfying the following conditions: 
\begin{enumerate}[noitemsep]
\item  The measure $m_x$ depends measurably on the point $x\in X$;
\item  Each measure $m_x$ has a finite first moment, i.e.,~for some (hence any) $o \in X$ and for any $x \in X$, $\int d(o,y)dm_x(y)$ is finite.
\end{enumerate}  

The data $(m_x)_{x\in X}$ allow us to define a notion of curvature as follows: 

\begin{definition}[{\cite[Definition 3]{Ollivier2007}}] 
\label{def:Ollivier}
Let $(X, d_X)$ be a metric space with a random walk $m$. 
Let $x,y \in X$ be two distinct points. 
The \emph{Ollivier-Ricci curvature} (briefly, \emph{OR curvature},  originally called the coarse Ricci curvature) of $(X,d_X,m)$ at $(x,y)$ is 
\begin{equation}\label{eq:OR}
\OR(x,y):=1-\frac{W_1(m_x ,m_y)}{d_X(x,y)}.
\end{equation}
\end{definition}

Similar to the classical Ricci case, this curvature will be positive or negative depending on whether the measures $m_x$ and $m_y$ are closer or further apart than the points $x$ and $y$ themselves. See~\cref{fig:OR} for a discussion on spaces with positive, zero or negative Ricci curvature; and their analogue graphs where curvature is measured using Ollivier-Ricci  curvature.

Consider a metric measure space $(X,d_X,\mu)$ where balls in $X$ have finite measure, and $\mu$ is fully supported on $X$. 
Choose some $\varepsilon>0$ and let $B(x,\varepsilon)$ denote the open $\varepsilon$-ball of $x\in X$. 
The \emph{$\varepsilon$-step random walk} on $X$, starting from a point $x$,  consists in randomly jumping in the ball of radius $\varepsilon$ around $x$,
with probability proportional to $\mu$, i.e.,~it is defined to be $m_x=\frac{\mu|_{B(x,\varepsilon)}}{\mu(B(x,\varepsilon))}$. 
In the case of unweighted graphs, $\varepsilon=1$ is a natural choice. 

Explicitly, for an unweighted graph $G$, the neighbourhoods of vertices $x$, $y$ adjacent to an edge $e = (x, y)$ are endowed with a uniform probability measure, i.e., for $z$ adjacent  to $x$ (denoted $z\sim x$),
\[\mu_x(z) := \frac{1}{\deg_G(x)}.\] 

In \cite{LinLuYau2011}, still focusing on the case of graphs, the authors proposed not to use a uniform measure but rather the following one: for any $t\in [0,1]$ and $x\in V$:  
\begin{equation}\label{eq:param_measure}
m_x^t(z):=\begin{cases}
    t, &\mbox{if $z=x$};\\
    \frac{1-t}{\deg_G(z)}, &\mbox{if $z\sim x$};\\
    0, &\mbox{otherwise.}
\end{cases}
\end{equation}
Here, the parameter $t$  controls whether the random
walk is likely to revisit a node. The original definition by Ollivier can be retrieved with $t={0}$.

The parametric definition for the measure can be further generalized as follows for a weighted graph (cf. \cite{BaiHuangLu2020}): for any $t\in [0,1]$ and $x\in V$,
\begin{equation}\label{eq:param_measure_weighted}
m_x^t(z):=\begin{cases}
    t, &\mbox{if $z=x$};\\
    \frac{(1-t)w_{xz}}{\sum_{z'\sim x} w_{xz'}}, &\mbox{if $z\sim x$};\\
    0, &\mbox{otherwise.}
\end{cases}
\end{equation} 
Here, $w_{xz}$ denotes the weight of the edge between vertices $x$ and $z$.

Another way of defining the measures for weighted graphs is by Ni et al.~\cite{NiLinLuo2019}: for any $t\in [0,1]$, $p\ge 0$ and $x\in V$, define the probability measure, 
\begin{equation}
\label{eqn:mx}
m_x^{t,p}(z):=\begin{cases}
    t, &\mbox{if $z=x$};\\
    \frac{1-t}{C_x}e^{-d_G(x,z)^p}, &\mbox{if $z\sim x$};\\
    0, &\mbox{otherwise.}
\end{cases}
\end{equation}
Here, $C_x=\sum_{z\sim x}e^{-d_G(x,z)^p}$ and $d_G$ is the graph distance. 
For any $x\sim z$, $d_G(x,z)$ is the weight between $x$ and $z$, with a default value of $1$ when the graph is unweighted.  
Note that when $t = 0$ and $p = 0$, or if the graph is unweighted, this measure reduces to the uniform measure. 
In this context, $t$ determines the probability of staying at $x$, and $p$ controls the extent to which the neighbour $z$ of $x$ is discounted based on the weight $d_G(x,z)$. %
The idea of the Ollivier-Ricci curvature definition compared to the Forman-Ricci curvature is illustrated in \cite[Figure 5]{TianLubbertsWeber2023}.
In the particular case when we use the probability measure following  \cref{eq:param_measure}, which depends on a parameter $t\in [0,1]$, we can consider the asymptotic curvature for $t$ tending to 1.

\begin{definition}[{\cite[Section 2]{LinLuYau2011}}]
\label{def:ORc-param}
    The $t$-Ollivier-Ricci curvature is
    \begin{equation}
    \label{eq:alpha-OR}
        \kappa_t(x,y):=1-\frac{W_1(m_x^t ,m_y^t)}{d_G(x,y)}. 
    \end{equation}
The \emph{LinLuYau-Ricci curvature}\footnote{We use this notion to differentiate it from the classical Ricci curvature.} at $(x,y)$ is
    \begin{equation}
    \label{eq:asymptotic-OR}
        \kappa^{{LLY}}(x,y):=\lim_{t\to 1}\frac{\kappa_t(x,y)}{1-t}.
    \end{equation}
\end{definition}
Following \cref{eq:asymptotic-OR}, we may consider a filtration of a graph using $t$ as a filtration parameter, giving rise to multiscale cores of a graph. 
This consideration is left for future discussions.  

\begin{remark}
We review examples of LinLuYau-Ricci curvatures and Forman-Ricci curvatures for various simple graphs. 
For the computation of LinLuYau-Ricci curvatures, see \cite[page 610]{LinLuYau2011} for details. For Forman-Ricci curvatures, the computation can be done using the formula $F(e)=4-\deg(u)-\deg(v)$, where $e$ is an edge connecting vertices $u$ and $v$. 
\begin{itemize}
    \item The complete graph $K_n$ has a constant LinLuYau-Ricci curvature on each edge equal to  $n/(n - 1)$. This is the only graph with a constant Ricci curvature greater than 1. 
    In $K_n$, the Forman-Ricci curvature is also a constant on each edge with a value of $6-2n.$
    \item The cycle $C_n$ for $n \ge 6$ has a constant LinLuYau-Ricci curvature on each edge equal to 0. For small cycles $C_3$, $C_4$, and $C_5$, we have constant LinLuYau-Ricci curvature on each edge equal to $3/2$, $1$, and $1/2$, respectively. 
    In $C_n$ for any $n\geq 3$, the Forman-Ricci curvature is a constant on each edge with a value of $0.$
    \item The hypercube $Q^n$ has a constant LinLuYau-Ricci curvature $2/n$ and a constant Forman-Ricci curvature $4-2n$.
\end{itemize}
\end{remark}

To deal with the lack of a resolution parameter of Ollivier-Ricci curvature and to reveal multiscale structures in real-world networks, Gosztolai and Arnaudon \cite{GosztolaiArnaudon2021} further modified  the definition of the probability measure. 
They introduced the idea of initiating a diffusion process at each node $i$ to generate a set of measures $p_i(\tau)$, where $\tau$ represents the time parameter of the diffusion process. 
The \emph{dynamic Ollivier-Ricci curvature} of an edge is then defined as the distance between the pair of measures started at its endpoints, normalized by the weight of the edge:
\begin{equation}
\label{eqDynOR}
\kappa_\tau(u,v):=1-\frac{W_1(p_u(\tau),p_v(\tau))}{w_{uv}}, 
\end{equation}
whenever $e =(u,v)$ is an edge, and $0$ otherwise. 

Following \cite{TianLubbertsWeber2023}, for a graph $G=(V,E)$, we can further define the Ollivier-Ricci curvature for vertices (rather than edges) with
respect to the curvature of its adjacent edges. Formally, if for a vertex $x\in V$ the set  $E_x := \{e \in E : x \in  e\}$ contains the edges adjacent to $x$, then the \emph{Ollivier-Ricci curvature at a vertex} $x$  is given by $\kappa(x) = \sum_{e\in E_x} \kappa(e)$. 

\para{Definitions for hypergraphs.}
To generalize the Ollivier-Ricci curvature to hypergraphs, Eidi et al.~\cite{EidiJost2020}  observed that \cref{def:Ollivier} can also be expressed by the following formula for an edge $e = (v,w)$ of a graph:
\begin{equation}
\label{eq:transport}
    \kappa{(e)} = \mu_0 - \mu_2 - 2\mu_3, 
\end{equation}
where $\mu_0$ represents the amount of mass that remains unmoved in an optimal transport plan, corresponding to the stable mass in directed 3-cycles $(u \to v \to w \to u)$ for any vertex $u$ adjacent to $v$ and $w$.  
The terms $\mu_2$ and $\mu_3$ are the amount of mass that should be moved with distance 2 and 3, respectively.

The above formula \cref{eq:transport} is obtained as follows:
consider the edge $e = (v,w)$  and let $e_v = (v, v_1)$ and $e_w = (w,w_1)$ be
edges emanating from $v$ and $w$, respectively. 
Then define
their distance w.r.t.~$e$ as $d_e(e_v,e_w):=d(v_1,w_1)$
where $d(v_1, w_1)$ denotes the distance between $v_1$ and $w_1$ in the graph, specifically the minimum number of edges that must be traversed to travel from $v_1$ to $w_1$ (by construction, not greater than 3:  distance 0,  when $v_1,w_1$ participates in a triangle, distance 1 for a square, and so on up to 3). 
Let $E_v$ be the set of edges that have
$v$ as a vertex, and let $|E_v|$ be its cardinality. 
Define a probability measure $m_v$ on the set of all edges $E$ of the graph by assigning each edge $e_v \in E_v$ a weight of $1/|E_v|$, and assigning a weight of 0 to all edges not in $E_v$. 
For an edge $e = (v,w)$, from $d(v,w)=1$ it follows that the Ollivier–Ricci curvature of \cref{def:Ollivier} is equal to 
\begin{equation}\label{eq:kappa}
 \kappa(e)=1-W_1(m_v,m_w),    
\end{equation}
where $W_1$ is the 1-Wasserstein distance between $m_v$ and $m_w$:
\begin{equation}\label{eq:wasserstein}
W_1(m_v,m_w)=\inf_{p\in\Pi(m_v,m_w)}\sum_{(e_1,e_2)\in E\times E}d_e(e_1,e_2)p(e_1,e_2), 
\end{equation}
where $\Pi(m_v,m_w)$ is the set of measures on $E \times E$ that project to $m_v$ and $m_w$, respectively. 
The key idea is to optimally arrange the two collections 
$E_v$, $E_w$ of edges sharing one of their endpoints with $e$, so that the average distances between the paired edges are minimized. 
Here, the sets $E_v$ and $E_w$ both include the edge $e = (v,w)$ that we are
evaluating. 
Let $\mu_i$ ($0\le i\le 3$) be the fraction of edges in $E_v$ that are moved at distance $i$ in some optimal transport plan.
According to \cref{eq:wasserstein}, $W_1(m_u, m_v)=\mu_1+2\mu_2+3\mu_3$. 
Observing that $\sum_{i=0}^3\mu_i=1$, \cref{eq:kappa} then yields \cref{eq:transport}.

The above formula \cref{eq:transport} is taken as the definition for the curvature of a hyperedge in a directed hypergraph in \cite{EidiJost2020}. 
Precisely, if we denote $u\to e_i$ when there exists a hyperedge $e = (e_k, e_i)$ such that $u \in e_k$, and  similarly, $e_j \to v$ when there is a hyperedge $e = (e_j, e_k)$ such that $v \in e_k$, then for a given hyperedge $e = (e_i, e_j)$, we can define two sets: ${\mathcal M} = \{u \colon u \to e_i\}$ referred to as the \emph{masses} and ${\mathcal H} = \{v \colon e_j \to v\}$, referred to as the \emph{holes}.  
A probability measure is assigned to each set, denoted by $\mu_{\mathcal M}$ and $\mu_{\mathcal H}$, respectively. 
For $u \in {\mathcal M}$ and $v \in{\mathcal H}$, the distance $d(u,v)$ between each mass $\mu_{\mathcal M}(u)$ and each hole $\mu_{\mathcal H}(v)$ of a hyperedge $e$ is defined as the minimum number of directed hyperedges connecting them. 
This distance is at most 3: it is exactly 3 if the shortest path from $u$ to $v$ requires passing through $e$; it is 0 when $u = v$, meaning that $u$ is both a mass and a hole of $e$.

A \emph{transport plan} is a matrix $\mathcal E$ whose entries represent the amount of mass, out of $\mu_{\mathcal M}(u)$, to be moved from vertex $u$ to vertex $v$, denoted by ${\mathcal E} (u, v)$. 
Assume we are given an optimal transport plan ${\mathcal E}$, i.e.,~one that minimizes the quantity $\sum_{u\to e_i}\sum_{e_j\to v}d(u,v){\mathcal E}(u,v)$. 
Let $\mu_i$ denote the amount of mass that is moved at distance $i$, with  $i\in \{0, 1, 2, 3\}$. 
The Ollivier-Ricci curvature $\kappa$ of $e$ for a hypergraph is then defined as $\kappa(e) = \mu_0 - \mu_2 - 2\mu_3$. It is bounded above by $\kappa = 1$ (achieved when $\mu_0 = 1$, meaning each mass coincides with a hole of its same size) and below by $\kappa = -2$ (achieved when $\mu_3 = 1$, meaning each mass must be moved a distance of 3). 

Coupette et al.~\cite{CoupetteDalleigerRieck2023} presented a flexible and unifying framework, ORCHID, which generalizes the Ollivier–Ricci curvature to hypergraphs by extending both the measure $\mu$ and the distance metric in several ways. The distance metric is replaced by a function aggregating measures (\texttt{AGG}):
\begin{equation} 
\kappa^{AGG}(u,v):=1-\frac{\texttt{AGG}(\mu_u,\mu_v)}{d(u,v)}. 
\end{equation}
and they identified three ways to generalize the probability measure induced by the lazy random walk on graphs: the \emph{equal-nodes random walk} on the unweighted clique expansion of $H$ that results from picking a neighbour $v$ of $u$ uniformly at random; the \emph{equal-edges random walk} that follows a two-step random walk on the unweighted star expansion of $H$; and the \emph{weighted-edges random walk} that can be thought as a random walk on a weighted clique expansion of $H$.

Treating a hyperedge as a set of nodes, the \texttt{AGG} function can be taken to be:
$$\texttt{AGG}_\mathrm{A}(e):=\frac{2}{|e|(|e|-1)}\sum_{\{u,v\}\subseteq e}W_1(\mu_u,\mu_v),$$
where $W_1(\cdot,\cdot)$ corresponds to the Wasserstein distance; $\texttt{AGG}_\mathrm{A}$ measures the average amount of effort to transport the probability mass from one node to another in the edge $e$.  In general, the \texttt{AGG} function should reduce to the Ollivier–Ricci curvature on graphs, while remaining computationally tractable and invariant under permutations of node indices.

The curvature at node $u$, independent of the choice for \texttt{AGG}, can be obtained as the mean of all curvatures of edges containing $u$: 
$$\kappa^{{E}}(u):=\frac{1}{\mathrm{deg}(u)}\sum_{e\ni u}\kappa^{AGG}(e).$$
Since the hypergraph $H$ is connected, Coupette et al.~\cite{CoupetteDalleigerRieck2023} defined the curvature of an arbitrary subset of nodes $s\subseteq V$ to be:$$\kappa^{AGG}(s):=1-\frac{\texttt{AGG}(s)}{d(s)},$$
where $d(s):=\max \{d(u,v)|\{u,v\}\subseteq s\}$. 
Coupette et al.~\cite{CoupetteDalleigerRieck2023} observed that the above definition aligns with their notions of hypergraph curvature for $s \in E, d(s)=1$.

\para{Definitions for simplicial complexes.}
In \cite{Yamada2023}, Yamada generalized the definition of Ollivier-Ricci curvature to weighted simplicial complexes. 
Let $S$ be a weighted simplicial complex with weight function $w$. 
Define the degree of a $p$-simplex $\alpha$ to be
\[\deg(\alpha):=\sum_{\alpha<\gamma\in S_{p+1}}w_\gamma.\]

Two distinct $p$-simplices $\alpha,\beta$ in $S$ are \textit{connected} if they share a $(p+1)$-coface $\gamma$, denoted as $\alpha\stackrel{\gamma}{\sim}\beta.$ If $\gamma$ exists, it is unique.
A path from a $p$-simplex $\alpha$ to another $p$-simplex $\beta$ is a finite sequence of connected $p$-simplices $\alpha=\alpha_0,\dots,\alpha_m=\beta$.
We call $m$ the length of the path.
We define the distance between two $p$-simplices $\alpha$ and $\beta$ to be the length of the shortest path from $\alpha$ to $\beta$, and denote it by $d_S(\alpha,\beta)$.

\begin{definition}[{\cite[Definitions 2.8 and 2.9]{Yamada2023}}]
    Let $S$ be a weighted simplicial complex with weight function $w$. 
    For a $p$-simplex $\alpha$ in $S$, define the probability measure
    \[m_\alpha(\alpha'):= \begin{cases}
        \frac{w_\gamma}{(p+1)\deg(\alpha)}, &\mbox{if $\alpha\stackrel{\gamma}{\sim}\alpha'$ for some $\gamma\in S_{p+1}$;}\\
        0, &\mbox{otherwise.}
    \end{cases}\]
    For any two distinct $p$-simplex $\alpha$ and $\beta$ in $S$, the Yamada-Ricci curvature at $(\alpha,\beta)$ is defined as
    \[\OR(\alpha,\beta):=1-\frac{W_1(m_\alpha,m_{\beta})}{d_S(\alpha,\beta)}.\]
\end{definition}

\para{Properties.}
When $X$ is a Riemannian manifold and $\varepsilon$ is small, the Ollivier-Ricci curvature derived from the $\varepsilon$-step random walks captures the Ricci curvature in the following manner. 
Consider $X$ a smooth, complete, $N$-dimensional Riemannian manifold. 
Let $x\in X$ and $v$ be a unit tangent vector at $x$. 
Let $y$ be a point on the geodesic that starts from $v$. 
As stated in \cite[Example 7]{Ollivier2009},
\[\OR(x,y)=\frac{\varepsilon^2\Ric(v,v)}{2(N+2)}+O(\varepsilon^3+\varepsilon^2 d(x,y)),\]
where $\Ric$ denotes the Ricci curvature and $O(\cdot)$ represents the big $O$ notation. 

The Ollivier-Ricci curvature bounds eigenvalues of the Laplacian of the graph (cf. \cite{Ollivier2009}, see also \cite{BauerJostLiu2012}), in the sense that if $k$ is a lower bound for the Ollivier-Ricci curvature of edges (i.e.,~$\kappa(e) \ge k$, for all edges $e$), then the eigenvalues $\lambda_i$ satisfy
\begin{equation}\label{eq:bound-gap}
     k \le \lambda_1 \le \cdots \le \lambda_{N-1} \le 2 - k,
\end{equation}
where $\lambda_1$ is the first non-zero eigenvalue.
It follows that the Ollivier-Ricci curvature provides formal bounds on the local clustering coefficient, which measures the connectivity of the network around a specific node. 
Thus, in the case of graphs, the Ollivier-Ricci curvature controls the amount of overlap between neighbourhoods of adjacent vertices. 

In other words, the Ollivier-Ricci curvature gives us information about how close we are to increasing the number of connected components in a graph: the smaller the eigenvalue gap, the ``thinner'' the bottleneck among two large components.
This observation motivates the application of the Ollivier-Ricci curvature to community detection, finding bottlenecks, and clustering (cf.~\cref{sec:graph-core-OR}). In \cite{EidiMukherjee2023}, Eidi and Mukherjee focused on this bridge between topology and geometry via random walks and the Laplacian. In particular, they considered random walks on $p$-simplices where $p\ge 1$ with the goal of connecting the spectral gap of the Laplacian in degree $p$ to the
$p$-dimensional homology. 
The authors then speculated about a connection with the Ollivier-Ricci curvature of simplicial complexes, which could potentially pave the way for using this curvature to detect how close a simplicial complex is to exhibiting higher-dimensional holes. 
  
According to Gosztolai and Arnaudon~\cite{GosztolaiArnaudon2021}, Ollivier-Ricci curvature measures how much a graph differs from being locally grid-like, analogous to the concept of being \emph{flat} in continuous spaces. 
The flatness of a network can be understood in terms of its local connectivity: the distance of a pair of nodes is the same as the average distance of their neighbourhoods. Thus, an edge with positive (or negative) Ollivier-Ricci curvature suggests that it is located in a part of the graph that is more (or less) connected than a grid.

Regarding the dynamic Ollivier-Ricci curvature in \cref{eqDynOR}, initially, as $\tau$ approaches zero---when all nodes support disjoint point masses and the diffusion has not yet mixed---the dynamic curvature tends toward $0$. 
At the other extreme, as $\tau$ approaches infinity, i.e.,~as the diffusion reaches a stationary state, the dynamic curvature tends to 1. 
At intermediate scales, the curvature can vary between 1 and some finite negative number, depending on the structure of the graph. 
According to \cite{GosztolaiArnaudon2021}, as the curvature of an edge evolves, the scale at which it approaches unity reflects how easily information can be propagated between clusters.

Interestingly, for applications of the Ollivier-Ricci curvature of hypergraphs to sparsification, removing vertices from a hyperedge  has the following effect on the curvature:
\begin{proposition}[{\cite[Proposition 2.4]{EidiJost2020}}]
Given a hyperedge $e \colon A = \{x_1, \ldots , x_n\} \to B = \{y_1, \ldots , y_m\}$, by removing $l\le n$ vertices from the set $A$ and $l'\le m$ vertices from B,  the following relation holds 
between the curvature of the resulting hyperedge $e'$ and $e$:
\[|\kappa(e)-\kappa(e')|\le \min\{3, 3\left(l/n,l'/m\right)\}.\]
\end{proposition}

However, in the graph setting, there are lower and upper bounds for the Ollivier-Ricci curvature of an edge that enjoy a combinatorial formula \cite{TianLubbertsWeber2023a}: a lower bound for the Ollivier-Ricci curvature $\kappa(e)$ of the edge $e$ with vertices $u$ and $v$ is given by
 \[\kappa^{low}(e):=-2\left(1-\frac{1}{\deg_G(u)}-\frac{1} {\deg_G(v)}\right),\]
 whereas an upper bound is given by
 \[\kappa^{up}(e):=\frac{\Delta(e)}{\max\{\deg_G(u),\deg_G(v)\}}.\]
 where $\Delta(e)$ denotes the number of triangles in $G$ that include the edge $e$.
 
Given these bounds,  the arithmetic mean of the lower and upper bounds yields  an approximation of  the Ollivier-Ricci curvature computable with a combinatorial formula:
\begin{equation}
\label{eq:ORapprox}
    \hat\kappa(e)=\frac{1}{2}\left(\kappa^{low}(e)+\kappa^{up}(e)\right).
\end{equation}

\subsection{Resistance Curvature}
\label{sec:resistance-curvature}

We now discuss a definition of curvature that is derived from the physical property of electric circuits called \emph{effective resistance}.
The effective resistance between two vertices in a weighted graph $G = (V,E, w)$ can be computed via the graph Laplacian \cite{KleinRandic1993}, which we review here. Let \(n=|V|\) and \(m = |E|\). Using the notation of Spielman and Srivastava~\cite{SpielmanSrivastava2011}, we first define $W$ to be an \(m\times m\) diagonal  matrix containing the weights $w$ of $G$. We then define the \emph{signed edge-vertex incidence matrix} $B \in \Rspace^{m \times n}$ as 
\begin{align}
B(v,e)  = \begin{cases}
0 & \textrm{if vertex $v$ is not on the boundary of edge $e$,}\\
1 & \textrm{if vertex $v$ is edge $e$'s head,} \\
-1 & \textrm{if vertex $v$ is edge $e$'s tail}. 
\end{cases}
\end{align} 
From here, we can define the \emph{graph Laplacian} $L \in \Rspace^{n \times n}$ as
\begin{align} 
L = B^\top W B, 
\end{align}
noting that $L$ is symmetric.

Letting $L^+$ denote the Moore-Penrose pseudoinverse of $L$, i.e.,~the unique matrix $L^+$ such that $LL^+ = L^+ L = \proj(\ker(L)^\perp)$, Spielman and Srivastava~\cite{SpielmanSrivastava2011} defined the \emph{effective resistance} of an edge $e$ to be the entry on the diagonal corresponding to $e$ of the matrix $R \in \Rspace^{m \times m}$ defined by   
\begin{align}
\label{eq:effective_resistances}
R := B L^{+} B^\top = B (B^\top W B)^+ B^\top.
\end{align} 
That is, the \emph{effective resistance} of and edge $e$ is given by:
$$r_e := R(e,e).$$

The pseudoinverse $L^+$ can also be used to define a notion of effective resistance between pairs of vertices~\cite{KleinRandic1993}:
\[r_{v_iv_j} = (\mathbf{e}_i-\mathbf{e}_j)^\top L^+(\mathbf{e}_i-\mathbf{e}_j)\]
where $\mathbf{e}_i \in \Rspace^n$ is a column vector of length $n =|V|$ with a one  in the $i$th spot and zeros everywhere else. 
Alternatively, if we let $\Gamma = L+\frac{J}{n}$, where $J$ is the $n\times n$ matrix of all ones, we have that the effective resistance between two vertices $v_i$ and $v_j$ is given by
\[r_{v_iv_j} = (\mathbf{e}_i-\mathbf{e}_j)^\top\Gamma^{-1}(\mathbf{e}_i-\mathbf{e}_j),\]
since $\Gamma$ in this case is invertible. 
We denote by \(\Omega\) the matrix whose \((i,j)\)-th entry is \(r_{v_i v_j}\); this matrix is known as the \emph{effective resistance matrix}. 

Devriendt and Lambiotte~\cite{DevriendtLambiotte2022} introduced discrete curvatures on both vertices and edges of weighted graphs based on effective resistance. 
They called them \emph{nodal} and \emph{link resistance curvature}, respectively. 
Given a weighted graph $G=(V,E,w)$, we denote these curvatures as $k^{up}: V \to \Rspace$ and $k^{down}: E \to \Rspace$, respectively. 
\begin{definition}(Nodal Resistance Curvature~\cite{DevriendtLambiotte2022}) 
\label{def:nodal_resistance_curvature}
The \emph{nodal resistance curvature} of a vertex $v_i \in V$ in a graph $G$ is defined as
\begin{align}
\label{eq:nodal-rc}
    k^{up}(v_i):=1-\frac{1}{2} \sum_{j \sim i} w_{v_iv_j}r_{v_iv_j}.
\end{align}
Here, $r_{v_iv_j}$ denotes the effective resistance between vertices $v_i$ and $v_j$, and $w_{v_iv_j}$ the weight of the edge $(v_i,v_j)$, where $j\sim i$ is the set such that there is an edge $e \in E$ connecting $v_i$ and $v_j$.
\end{definition} 
Devriendt and Lambiotte~\cite{DevriendtLambiotte2022} called the quantity $w_{v_iv_j}r_{v_iv_j}$ the \emph{relative resistance} of the edge $e = (v_i,v_j)$ as it quantifies how important the edge $e$ is for the connectivity of $G$. 

We note that Devriendt et al.~\cite{DevriendtOttoliniSteinerberger2024} introduce a slightly different definition of resistance curvature, namely
\begin{align}
k^{up} = \Omega^{-1}\mathbf{1}.  
\end{align}
Observe that the expression in \cref{eq:nodal-rc} can also be written as
\begin{align}
 k^{up} = \frac{{\Omega}^{-1}\mathbf{1}}{\langle \mathbf{1},\Omega^{-1}\mathbf{1}\rangle}.   
\end{align}
It is evident the two definitions are equal up to a factor of $\langle \mathbf{1},\Omega^{-1}\mathbf{1}\rangle$.

As with the curvatures mentioned in previous sections, there is also a definition of resistance curvature on edges: 
\begin{definition}(Link Resistance Curvature~\cite{DevriendtLambiotte2022}) 
\label{def:link_resistance_curvature}
The link resistance curvature of an edge $e=(v_1,v_2) \in E$ in a graph $G$ is defined as 
\begin{align}
    k^{down}(e)=\frac{2(k^{up}(v_1)+k^{up}(v_2))}{r_{v_1v_2}}.
\end{align}
\end{definition}
Intuitively speaking, the resistance curvature of a vertex is derived from the effective resistance of its adjacent edges, whereas the resistance curvature of an edge comes from the resistance curvature of its vertices. Let $\kappa$ denote the Ollivier-Ricci curvature and $\mathcal{F}$ denote the Forman-Ricci curvature, it was shown in \cite{DevriendtLambiotte2022} that for an edge $e = (v_i, v_j)$,
\[\kappa(e) \geq k^{down}(e) \geq \frac{\mathcal{F}(e)}{\omega_{e}}.\]
We believe it is easy to generalize resistance curvature to hypergraphs and simplicial complexes, which is left for future work.

\begin{table}[!t]
\caption{A summary of recent works in finding cores of graphs, as surveyed in \cref{sec:graph-cores}. 
Ref.: references. V: vertices (nodes). E: edges. FRC: Forman--Ricci curvature. A-FRC: augmented FRC. ORC: Ollivier--Ricci curvature. ER: effective resistance. NMI: normalized mutual information. ARI: adjusted Rand index. MAE: mean absolute error. 
``Threshold": a subsampling method that removes vertices or edges above/below a chosen value.
``Probability'': a subsampling method that removes vertices or edges based on a probability distribution.
``Spectral'': properties related to the Laplacian matrix of a graph.}
\label{table:graph-core-summary}
\centering
\begin{tabular}{p{1cm}p{1.5cm}p{2.5cm}p{2.5cm}p{3.5cm}}
\hline\noalign{\smallskip}
\textbf{Ref.} & \textbf{Graph} & \textbf{Method} & \textbf{Criteria} & \textbf{Property Preserved} \\
\noalign{\smallskip}\svhline\noalign{\smallskip}
\cite{WeberJostSaucan2018}         & weighted    & threshold, E             & FRC                & backbone \\
\cite{BarkanassJostSaucan2022}     & weighted    & threshold, V or E        & FRC                & backbone \\
\cite{SalamatianAndersonMatthews2022} & weighted & threshold, E             & FRC                & backbone \\
\cite{ZhangSongTao2024}            & unweighted  & threshold, E             & Balanced FRC       & clusters \\
\cite{KimJeongLim2022}             & unweighted  & threshold, E             & FRC                & clusters (modularity, NMI) \\
\cite{FesserIvanezDevriendt2024}   & unweighted  & threshold, E             & A-FRC              & clusters (accuracy) \\
\cite{TianLubbertsWeber2023}       & weighted    & threshold, E             & FRC, ORC, A-FRC    & single- \& mixed-membership clusters (NMI) \\
\cite{ParkLi2024}                  & unweighted  & threshold, E             & lower FRC          & clusters \\
\noalign{\smallskip}\hline\noalign{\smallskip}
\cite{SiaJonckheereBogdan2019}     & weighted    & threshold, E             & ORC                & clusters (accuracy) \\
\cite{NiLinLuo2019}                & weighted    & threshold, E             & ORC                & clusters (modularity, ARI) \\
\cite{WuChengCai2023a}             & weighted    & threshold, V and E       & ORC                & clusters (MAE) \\
\cite{GosztolaiArnaudon2021}       & weighted    & N/A                      & dynamical ORC      & multiscale clusters (geometric modularity) \\
\cite{TianLubbertsWeber2022}       & weighted    & threshold, E             & ORC                & clusters (NMI) \\
\noalign{\smallskip}\hline\noalign{\smallskip}
\cite{SpielmanSrivastava2011}      & weighted    & probability, E           & ER                 & spectral \\
\noalign{\smallskip}\hline
\end{tabular}
\end{table}

\section{Finding Cores of Graphs}
\label{sec:graph-cores}
In this section, we survey recent approaches in finding cores of graphs using geometric or topological means. We begin with the applications of discrete curvatures, in particular, Forman-Ricci Curvature (\cref{sec:graph-core-FR}) and Ollivier-Ricci Curvature (\cref{sec:graph-core-OR}), followed by a discussion of graph sparsification based on effective resistance (\cref{sec:graph-core-ER}). 
References surveyed in this paper are classified in \cref{table:graph-core-summary}.

\subsection{Forman-Ricci Curvature}
\label{sec:graph-core-FR} 

A number of Forman-Ricci curvature-based graph sparsification methods have been proposed to preserve the \emph{backbone} (e.g.,~\cite{WeberJostSaucan2018, BarkanassJostSaucan2022, SalamatianAndersonMatthews2022, ZhangSongTao2024}). 
Weber et al.~\cite{WeberJostSaucan2018} defined the backbone of a network (graph) as a subnetwork (subgraph) that captures important nodes (called \emph{hubs}) and edges (called \emph{bridges}). 
Hubs are nodes with high degree and high betweenness centrality, whereas bridges are edges governing the mesoscale structure of the graph (such as those forming long-range connections)~\cite{BarkanassJostSaucan2022}.  
According to Barkanass et al.~\cite{BarkanassJostSaucan2022}, a backbone is considered to be \emph{structure-preserving} if it preserves the structural features such as the node degree distribution or community structure of the original graph. 
The notions of a \emph{network backbone} and the \emph{core of a network} are sometimes used interchangeably~\cite{Newman2018}. 
In this survey, however, we define the \emph{core of a higher graph} more generally as a minimalist representation that retains its geometric or topological information. 

On the other hand, graph clustering or community detection can also be regarded as core finding by grouping (thus preserving) the highly interconnected regions of a graph, while discarding its  small-scale information. 
To that end, the Forman-Ricci curvature has also been used as a criterion to guide graph clustering (e.g.,~\cite{FesserIvanezDevriendt2024,TianLubbertsWeber2023,ParkLi2024}). 

\subsubsection{Curvature-Based Sampling for Preserving Backbones}
\label{sec:sparsification-fr}

Using the Forman-Ricci curvature, Weber et al.~\cite{WeberJostSaucan2018} proposed a sampling method to reduce the size of a given graph while preserving its backbone. 
The key idea is that ``high absolute curvature is strongly related to the structural importance of an edge.''~\cite{WeberJostSaucan2018}. 
Their sampling procedure applies thresholding with respect to the Forman-Ricci curvature, where edges with high absolute curvature are selected and edges with low absolute curvature are eliminated. 
They experimentally showcased that some desirable structure of the graph (such as clusters or communities) are preserved after the sampling. As an example, they showed that in the weighted graph of co-occurrences in Les Mis\'erables, the backbone graph maintains relationships between Valjean and other characters central to the storyline, and preserves a cluster of revolutionaries around Gavroche. 

Barkanass et al.~\cite{BarkanassJostSaucan2022} further extended the above sampling procedure by considering Forman-Ricci curvatures for graphs and 2-dimensional simplicial complexes, as well as the Haantjes-Ricci curvature, i.e., a form of curvature in general metric spaces. 
They showed that sampling by discrete versions of Ricci curvature is mathematically justified by the relation between the Ricci curvature of a manifold and the Forman curvatures of its discretization. 
They further studied coarse embeddings of a network into a metric space using kernels based on Forman-Ricci curvature.  

For real-world applications, Salamatian et al.~\cite{SalamatianAndersonMatthews2022} discussed challenges in accurately assessing network connectivity within large cloud provider infrastructures, due to potential manipulation of probe packets. They demonstrated the success of the above sampling method based on Forman-Ricci curvature, using latency measurements from RIPE Atlas anchors and virtual machines in data centers of three major cloud providers.

\subsubsection{Graph Sparsification for Graph Learning}
\label{sec:graph-core-FR-graph-learning}
By employing Ricci curvature-based sampling, Zhang et al.~\cite{ZhangSongTao2024} proposed a graph sparsification method based on the \emph{balanced Forman-Ricci curvature}~\cite{ToppingDiGiovanniChamberlain2021}---a modification of the Forman-Ricci curvature---to study over-squashing in deep networks. 
The balanced Forman-Ricci curvature takes into account the presence of triangles and cycles of length four. 
To perform graph representation learning, Zhang et al.~stored the topological information in the form of a sparsified computation subgraph around each node, where such a sparsified subgraph is obtained by iteratively sampling the edges with the highest balanced Forman-Ricci curvature in neighbourhoods (of increasing radius) around the node.

For theoretical foundations of graph sparsification, Spielman and Teng~\cite{spielman2004nearly} introduced nearly-linear time algorithms for spectral sparsification, leveraging reduced edge sets to preserve graph Laplacians. 
Benczúr and Karger~\cite{benczur1996approximating} introduced randomized edge sampling to approximate minimum s-t cuts, establishing a practical sparsification technique. 
Later, Batson et al.~\cite{batson2012twice} proposed deterministic spectral sparsifiers and proved that every graph has a spectral sparsifier with a number of edges linear in its number of vertices.

Graph sparsification also appears frequently in graph neural networks (GNNs), where a number of approaches have been developed to leverage or learn sparse graph structures. 
Peng et al.~\cite{PengGurevinHuang2022} discussed two primary approaches to reducing the training and inference complexity of GNNs through sparsification: sparsifying either the input graph or the model itself. Specifically, they formulated model sparsification for GNNs using both the train-and-prune paradigm and sparse training. 
Chen et al.~introduced FastGCN~\cite{chen2018fastgcn} that addresses scalability by introducing an importance sampling technique that sparsifies the neighbourhood aggregation process, leading to efficient graph convolutional network (GCN) training. 
Jin et al.~\cite{JinMaLiu2020} proposed a method to jointly learn sparse graph structures and node features for robust GNN performance under noise. 
Franceschi et al.~\cite{franceschi2019learning} treated the graph structure as a discrete variable in an optimization setting, enabling the discovery of sparse dependency  structures during training. See~\cite{ZhuXuZhang2022} for a survey of graph structure learning.

\subsubsection{Forman-Ricci Curvature in Community Detection}
\label{sec:graph-core-FR-community}

A number of applications within the field of community detection and graph clustering \cite{KimJeongLim2022, FesserIvanezDevriendt2024,TianLubbertsWeber2023,ParkLi2024} employ the Forman-Ricci curvature as a pruning or clustering criterion.  

Kim et al.~\cite{KimJeongLim2022} proposed \emph{link pruning} as a preprocessing step for community detection. 
Their pruning procedure eliminates edges below a certain threshold according to edge attributes such as the Forman-Ricci curvature. 
The threshold is chosen to preserve a proportion of the edges in the graph, in a similar manner to the sampling algorithms described in \cref{sec:sparsification-fr}. 
After preprocessing, Kim et al.~performed community detection using standard methods such as the Louvain method \cite{BlondelGuillaumeLambiotte2008}.  
They used modularity (for real-world social networks) and normalized mutual information (NMI) for synthetic networks with ground-truth communities to evaluate the performance of these community detection methods. 
Roughly speaking, modularity measures the density of connections within communities, and NMI describes the similarity between the predicted communities with the ground-truth. 
The authors demonstrated empirically that community detection with link pruning in most cases achieved better performance in terms of modularity and NMI, in comparison with traditional graph sparsification. 

As we will see in~\cref{sec:graph-core-OR}, community detection using the Ollivier-Ricci curvature is well-founded and has been proven to perform well in synthetic and real-world data. These results are based on the fact that positively curved edges, with respect to the Ollivier-Ricci curvature, tend to lie within communities, whereas negatively curved edges tend to serve as bridges between communities. The Forman-Ricci curvature behaves similarly. 
Although the Forman-Ricci curvature is advantageous for computational purposes compared to the Ollivier-Ricci curvature, it may be too simple in some cases to perform well in community detection.   

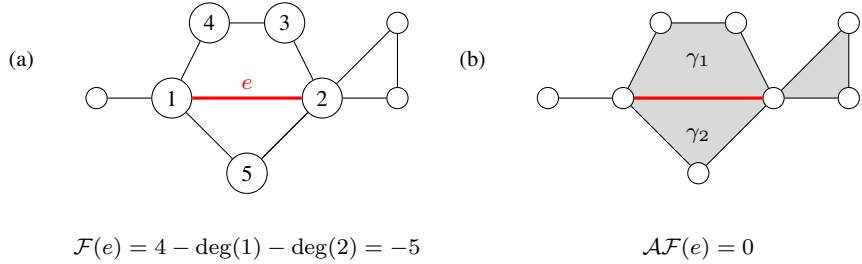
\begin{figure}[!htb]
    \centering
    \begin{tikzpicture}
    \begin{scope}[shift={(0,0)}]
        \node at (-1,1.5) {(a)};
        \foreach \i/\x/\y in {1/1/1, 2/3/1, 3/2.5/2, 4/1.5/2, 5/2/0} {
            \node[draw, circle] (v\i) at (\x,\y) {\i};
        }
        \foreach \i/\x/\y/\label in {6/0/1/, 7/4/1/, 8/4/2/} 
        {
            \node[draw, circle] (v\i) at (\x,\y) {\label};
        }
        \foreach \i/\j in {1/4, 1/6, 1/5, 2/5, 2/5, 3/4, 3/2, 2/7, 2/8, 7/8} 
        {
            \draw (v\i) -- (v\j);
        }
        \draw[red, very thick] (v2) -- node[above] {$e$} (v1);
        \node at (2,-1) {$\mathcal{F}(e) = 4 - \deg (1) - \deg(2) = -5$};
    \end{scope}

    \begin{scope}[shift={(6,0)}]
        \node at (-1,1.5) {(b)};

        \fill[gray!30] (1,1) -- (3,1) -- (2.5,2) -- (1.5,2) -- cycle;
        \fill[gray!30] (1,1) -- (3,1) -- (2,0) -- cycle;
        \fill[gray!30] (3,1) -- (4,1) -- (4,2) -- cycle;

        \foreach \i/\x/\y in {1/1/1, 2/3/1, 3/2.5/2, 4/1.5/2, 5/2/0, 6/0/1/, 7/4/1/, 8/4/2/} {
            \node[draw, circle, fill = white] (v\i) at (\x,\y) {};
        }
        \foreach \i/\j in {1/4, 1/6, 1/5, 2/5, 3/4, 3/2, 2/7, 2/8, 7/8} {
            \draw (v\i) -- (v\j);
        }

        \draw[red, very thick] (v2) -- node[above] {} (v1);

        \node at (2,-1) {$\mathcal{AF}(e) = 0$};

        \node at (2,1.5) {$\gamma_1$};
        \node at (2,0.5) {$\gamma_2$};

    \end{scope}
    \end{tikzpicture}
    \caption{Forman-Ricci curvature vs.~augmented Forman-Ricci curvature. Left: Forman-Ricci curvature of the edge \(e\). Right: augmented Forman-Ricci curvature of the same edge. We see that the edge is contained in the cycles $\gamma_1 = 1234$ and $\gamma_2 = 125$. Looking at \cref{eq:augmented_fr}, we then see that $\mathcal{AF}(e) = 2 + 2 - 3 -1 = 0.$  This image is a reproduction of Figure 1 from \cite{FesserIvanezDevriendt2024}, licensed under CC BY (Creative Commons Attribution) and generated with the assistance of ChatGPT.}
    \label{fig:fr_and_augmented_fr}
\end{figure}

Consequently, Fesser et al.~\cite{FesserIvanezDevriendt2024} proposed to use the \emph{augmented Forman-Ricci curvature}~\cite{ivanez2022} for community detection and explored  the connection between different augmentations and the expressive power and computational complexity for community detection. 
The augmented Forman-Ricci curvature is obtained from augmentations of the graphs. 
That is, by adding 2-dimensional cells as cofaces of all cocycles to obtain an unweighted cellular complex, and applying \cref{def:forman-ricci-sc} to the edges, we arrive at
\begin{equation}
\label{eq:augmented_fr}
    \mathcal{AF}(e) = 2 +\Gamma_{ee} - \sum_{e'\sim e}|\Gamma_{ee'} + \Gamma_{e'e} -1| - \sum_{e \not\sim e'} |\Gamma_{ee'}- \Gamma_{e'e}|,
\end{equation}
where $e\sim e'$ denotes two edges sharing a vertex and $e\not\sim e'$ two edges not sharing a vertex, and $\Gamma_{ee'}$ is the number of the cycles containing both \(e\) and \(e'\). See \cref{fig:fr_and_augmented_fr} for an example of the computation of this curvature. An alternative notion of augmented Forman-Ricci curvature can be obtained by restricting the length of the cycles considered in the definition of the augmentations of the graph:  
\[\mathcal{AF}_n(e) = \mathcal{AF}(e) \text{ taking into account all cycles of at most length }n.\]
This definition is supported by the idea that curvature should be a local construction. 
Fesser et al.~\cite{FesserIvanezDevriendt2024} mostly worked with $\mathcal{AF}_n$ with \(n \in \{3,4\}\) for community detection.

Several authors have studied the similarities between Ollivier-Ricci curvature and Forman-Ricci curvature \cite{TeeTrugenberger2021, JostMuch2021, SamalSreejithGu2018}, establishing empirical and theoretical equivalences and connections between these two notions. 
Fesser et al.~\cite{FesserIvanezDevriendt2024} performed an extensive correlation analysis of the augmented Forman-Ricci curvature against the Ollivier-Ricci and the Forman-Ricci curvatures. 
They applied augmented Forman-Ricci curvature to community detection and demonstrated performance improvement against the Ollivier-Ricci curvature-based method. 
In general, community detection methods using the augmented Forman-Ricci curvature are faster and less computationally expensive than those based on the Ollivier-Ricci curvature.

Tian et al.~\cite{TianLubbertsWeber2023} recently brought forward a unified framework for Ollivier-Ricci curvature-based community detection, including Ollivier and Forman-Ricci curvatures, in order to study their respective weaknesses and strengths. 
As part of this framework, they implemented a community detection algorithm based on the augmented Forman-Ricci curvature, with a weighting scheme for the edges (cf.,~\cite{FesserIvanezDevriendt2024}).
In practice, they implemented \(\mathcal{AF}_3\), that is, they considered only triangles and worked with simplicial complexes of dimension 2. 
They empirically showed that this reduction achieves a similar performance to the one with the full augmented Forman-Ricci curvature, while providing a significant speed-up in computations. In addition, they addressed a gap in the literature by incorporating mixed membership community detection to their analysis. In mixed-membership communities, one node might belong to more than one community, so that there might be overlaps between communities. Mixed membership is inspired by real-life situations. For instance, in a social network, a person might belong to more than one group of friends, such as school friends and university friends; or a protein might have several functional roles in a biochemical network (see \cite[Figure 1]{TianLubbertsWeber2023} for an illustrative example). In order to incorporate this framework with curvature-based community detection methods, Tian et al.~proposed to turn to the dual or line graph.

Park and Li~\cite{ParkLi2024} recently introduced the \emph{lower Ricci curvature} for community detection. This new notion of curvature is inspired by the balanced Forman-Ricci curvature (already mentioned in \cref{sec:graph-core-FR-graph-learning}), a curvature built upon the Forman-Ricci curvature to address its unboundedness and scale-dependent nature, in addition to its skewness toward negative values. As an extension of the Forman-Ricci curvature, the lower Ricci curvature has linear computational complexity. Park and Li also draw connections between the lower Ricci curvature and the Cheeger constant, which motivates the use of this curvature in community detection. 
In more detail, they proved that if there exists a positive constant \(\alpha >0\), that is, a lower bound of the lower Ricci curvature for all the edges in the graph, then the diameter of the graph is bounded from above by \(\alpha/2\) and the Cheeger inequality of the graph is bounded from below by the same quantity. This observation can be interpreted as follows: for larger values of \(\alpha\), the graph tends to be fully connected, hence having a smaller diameter and a larger Cheeger constant, implying more interconnectivity and less well-separated communities within the graph. Given these theoretical guarantees, Park and Li proposed a preprocessing algorithm that removes edges with lower Ricci curvature below a certain threshold, computed as the local minimum of a Gaussian mixture model fitted to the distribution of the lower Ricci curvatures. Using this preprocessing algorithm, they demonstrated improvement in efficiency and performance of a number of community detection methods across simulated and real-world data. 

\subsection{Ollivier-Ricci Curvature}
\label{sec:graph-core-OR}

The most popular application of the Ollivier-Ricci curvature for core finding is in clustering and community detection. 
Sia et al.~\cite{SiaJonckheereBogdan2019} explored community detection using the Ollivier-Ricci curvature (\cref{eq:OR}) on weighted graphs. 
They considered probability measures $m_x(v)=\frac{w_{xv}}{\sum_{v'\sim x} w_{xv'}}$ if $v$ is adjacent to $x$ and $m_x(v)=0$ otherwise. Here, $w_{xv}$  
denotes the weight of the edge between vertices $x$ and $v$. In their study, negatively curved edges act as \emph{bottlenecks} or \emph{bridges} linking different communities, whereas positively curved edges represent within-community connections.
Their algorithm removes edges with the most negative Ollivier-Ricci curvature, and then re-calculates the curvature for the affected edges. The process is repeated until all edges have non-negative curvature.
Compared with other community detection methods based on modularity or edge betweenness, this curvature-based method showed comparable or better accuracy.  

Ni et al.~\cite{NiLinLuo2019} introduced a family of probability distributions $m_x^{t,p}$ with $t\in [0,1]$ and $p\geq 0$; see \cref{eqn:mx}. 
In most computations, they set $t=0.5$ and $p=2$, and utilized the resulting random walk $\{m_x^{0.5,2}\}_x$ to define the Ollivier-Ricci curvature; see \cref{fig:OR} for an illustration of this curvature in the case of tree graphs, grid graphs and complete graphs.
They introduced a discrete Ricci flow for a weighted graph iteratively defined as $w_{xy}^{(i+1)}:=\left(1-\kappa^{(i)}(x,y)\right)d^{(i)}(x,y)$, with the superscript representing the iteration step, $d^{(0)}(x,y)=d(x,y)$, and $w_{xy}^{(0)}$ being the initial weight between $x$ and $y$. The algorithm runs by first applying the graph with 20 to 50 iterations of discrete Ricci flow processes, and then removing edges from highest curvature to lowest. 
The hierarchical community structure of the graph is revealed through careful adjustment of the cut-off threshold.
Experiments were conducted on both synthetic networks and real-world datasets with ground-truth communities.
Compared to several other community detection methods based on  modularity or edge betweenness, their algorithm demonstrates competitive or better performance, measured using modularity and adjusted Rand index (ARI), which measures the accuracy of clustering results against the ground-truth. 

Wu et al.~\cite{WuChengCai2023a} utilized the $t$-Ollivier-Ricci   
curvature of graphs (defined in \cref{eq:alpha-OR}) with $t=0.5$  to develop a graph sampling algorithm that reduces the number of nodes and edges while preserving the number of communities.
Their algorithm starts by randomly selecting an edge and calculating its Ollivier-Ricci curvature, then then expands the subgraph by choosing the next edge with the largest difference in Ollivier-Ricci curvature from the previously selected edge.
Their algorithm captures community structure by differentiating between within-community edges and between-community edges, which is supported by their finding that edges within a community have larger Ollivier-Ricci curvature than edges between that community and some other communities.
Compared with degree-based graph sampling methods via the mean absolute error (MAE), this curvature-based method 
is more effective at preserving minor communities.

Different from \cite{SiaJonckheereBogdan2019,NiLinLuo2019} where clusters are identified based on finding negatively curved edges between clusters, Gosztolai and Arnaudon~\cite{GosztolaiArnaudon2021} showed that along its evolution, the distribution of the dynamical Ollivier-Ricci curvature of edges exhibits gaps, i.e.,~differences in the relative magnitude of curvatures, which at characteristic timescales indicate bottleneck-edges that limit information spreading. 
These curvature gaps are claimed to be robust against large fluctuations in node degrees, preserving community structures up to the phase transition of detectability. 
As the diffusion process advances, these gaps capture progressively coarser community features; and they are used to uncover multiscale communities by analyzing deviations from constant curvature.

Tian et al.~\cite{TianLubbertsWeber2022} interpreted the graph clustering algorithms of Ni et al.~\cite{NiLinLuo2019} and Sia et al.~\cite{SiaJonckheereBogdan2019} as specific instances of a broader algorithm. 
These algorithms are based on geometric flows associated with discrete curvature \cite{Ollivier2009}, under which the edge weights of the graph evolve in a way that reveals its community structure.
Starting with a possibly unweighted graph, let the edge weights evolve under the curvature flow, up to a re-normalization at each iteration: calling $G^T=(V,E^T,w^T)$ the weighted graph at the iteration $T$, having fixed cut-off thresholds $h_T$, the new edge set is  
\[E^{T+1}:=\{e\in E^T: w^T(e)>h_T\},\] 
and the updated weights for any edge $e\in E^{T+1}$ with $e=\{u,v\}$,
\[w^{T+1}(e):=(1-\kappa(e))d_G(u,v),\]
where $\kappa$ is the discrete curvature of $G^T$ and similarly  $d_G$ is the graph distance of $G^T$.

The intuition is that, over time, the negative curvature of edges connecting different communities intensifies: edges with lower curvature experience a reduction in weight under the flow. 
On the other hand, edges with higher curvature get larger weights. Thus, interpreting weights as lengths of edges, the internal edges will contract faster over time, and the bridges will contract slower. The above algorithm is independent of the particular discrete curvature chosen. In any case, the success of the method depends crucially on identifying a good weight
threshold for cutting edges at each step. To speed up the algorithm in the case of the Ollivier-Ricci curvature, Tian et al.~used the combinatorial approximation of the Ollivier-Ricci curvature given in \cref{eq:ORapprox}.

\subsection{Graph Sparsification with Effective Resistance}
\label{sec:graph-core-ER}

Another avenue to find the core of a graph that preserves its spectral properties utilizes the notion of effective resistance (see \cref{eq:effective_resistances}), which is directly linked to the Laplacian of a graph.

Spielman and Srivastava~\cite{SpielmanSrivastava2011} introduced a nearly-linear time algorithm that, given a weighted graph \(G=(V, E, w)\), produces a sparsified graph \(H=(V, F, u)\) that approximately preserves the graph Laplacian of $G$. 
$H$ and $G$ share the same set of vertices, but $H$ contains fewer edges. 
Let $R(e,e)$ be the effective resistance of an edge $e$. The algorithm samples a predetermined number \(q\) of edges \(e\) from \(G\) with probability \(p(e)\) proportional to \(w(e) R(e,e)\) (the relative effective resistance, as named in \cite{DevriendtLambiotte2022}) and adds them to \(H\) with weight \(w(e,e) / (qp(e))\). The sampling is done independently with replacement, summing the weights if an edge is sampled more than once. Under such a construction, for \(0 <\varepsilon\leq 1\), the algorithm produces a sparsified graph $H$ after \(q\) sampling steps, such that the Laplacian of the sparsified graph \(L_H\) with (re-weighted) edges satisfies the following relation with respect to the Laplacian of the original graph \(L_G\):
\[(1-\varepsilon) \mathbf{x}^\top L_G \mathbf{x} \leq \mathbf{x}^\top L_H \mathbf{x} \leq (1+\varepsilon) \mathbf{x}^\top L_G \mathbf{x}, \quad \forall \mathbf{x} \in \Rspace^n.\] 
The Courant-Fischer theorem, provides a variational characterization of the eigenvalues \(\lambda_i\), \(1\leq i\leq n\) of the Laplacian \(L\) in terms of the following min-max expression: 
\[\lambda_i = \max_{\substack{V \subseteq \Rspace^n\\ \dim(V) = i}} \min_{\substack{\mathbf{x} \in V\\\mathbf{x} \neq \mathbf{0}}}\dfrac{\mathbf{x}^\top L \mathbf{x}}{\mathbf{x}^\top \mathbf{x}}.\]
The above expression allows us to conclude that
\[(1-\varepsilon) \, \lambda_i(G) \leq \lambda_i(H) \leq (1+\varepsilon)\, \lambda_i(G), \quad \forall\, 1\leq i \leq  n. \]
The strategy followed in this work was later extended to simplicial complexes by Osting et al.~\cite{OstingPalandeWang2020} (reviewed in~\cref{sec:sc-cores-spectral}).\\

Our survey of methods for finding graph cores shows that curvature-based approaches predominantly employ thresholding techniques to remove vertices or edges.
This preference likely stems from theoretical insights into discrete curvatures, which suggest that edges with relatively high or low absolute curvatures convey certain structural information, thereby justifying their removal based on threshold criteria.
On the other hand, the spectral method that leverages effective resistance employs a probability-based sampling strategy, which removes vertices or edges according to a predefined probability distribution. For future work, it would be interesting to extract graph cores by sampling with resistance curvature, investigate the properties preserved, and study the connection between effective resistance based sparsification and resistance curvature based sampling.   

A number of research works utilize probability-based sampling methods (often referred to as network percolation) to extract cores from graphs or to study the effect of random removal of nodes and edges. 
Since network percolation is a well-established field, we do not review relevant works in this paper, instead, we refer readers interested in this topic to a recent survey \cite{Sahimi2023Applications}. 
On the other hand, we believe that hypergraph percolation for core finding is a fairly new topic and deserves some attention, which is surveyed in~\cref{sec:hypergraph-cores}.

\begin{table}[!t]
\caption{A summary of recent works in finding cores of hypergraphs, as surveyed in \cref{sec:hypergraph-cores}. 
V: vertices (nodes). E: hyperedges. ER: effective resistance. ``robustness'': hypergraph robustness.}
\label{table:hypergraph-core-summary}
\centering
\begin{tabular}{p{1cm}p{2cm}p{2.7cm}p{2.3cm}p{3cm}}
\hline\noalign{\smallskip}
\textbf{Ref.} & \textbf{Hypergraph} & \textbf{Method} & \textbf{Criteria} & \textbf{Property Preserved} \\
\noalign{\smallskip}\svhline\noalign{\smallskip}
\cite{SunBianconi2021}            & undirected     & percolation, $V$ and $E$      & none                              & robustness                        \\
\cite{LeeGohLee2023}              & undirected     & percolation, $V$ and $E$      & degree                            & robustness, subhypergraph         \\
\cite{BianconiDor2024}            & undirected     & percolation                  & none                              & robustness, subhypergraph         \\
\noalign{\smallskip}\hline\noalign{\smallskip}
\cite{SomaYoshida2019}            & (un)directed   & sparsification               & Laplacian                         & spectral                          \\
\cite{KapralovKrauthgamerTardos2021}& (un)directed  & sparsification               & Laplacian                         & spectral                          \\
\cite{Lee2023}                    & undirected     & sparsification               & Laplacian                         & spectral                          \\
\cite{JambulapatiLiuSidford2023}  & undirected     & sparsification               & Laplacian                         & spectral                          \\
\cite{ChekuriXu2018}              & undirected     & cut sparsification           & max cut                           & global maximum cut                \\
\cite{ChenKhannaNagda2020}        & undirected     & cut sparsification           &  auxiliary graph \; \; \; edge strength     & spectral                          \\
\cite{AghdaeiFeng2022}            & weighted       & threshold, $V$ and $E$        & ER                                & clusters                          \\
\noalign{\smallskip}\hline\noalign{\smallskip}
\cite{ZhouHuangScholkopf2006}     & undirected     & normalized cut               & Laplacian                         & clusters                          \\
\cite{YangDengLu2021}             & weighted       & normalized cut               & intra-weights, $E$                & clusters                          \\
\cite{CoupetteDalleigerRieck2023} & undirected     & threshold, $E$               & ORC, edge and \; \; \; node features       & clusters                          \\
\noalign{\smallskip}\hline
\end{tabular}
\end{table}

\section{Finding Cores of Hypergraphs}
\label{sec:hypergraph-cores}

Numerous studies have focused on identifying the cores of hypergraphs. Techniques include finding the $(k,q)$-core of the hypergraph~\cite{LeeGohLee2023}, preserving the \emph{quadratic form of the Laplacian} of the resulting hypergraph~\cite{SomaYoshida2019}, and approximately preserving the cuts~\cite{ChekuriXu2018}; see~\cref{table:hypergraph-core-summary} for a summary. 
There is a plethora of clustering and community detection studies for hypergraphs (see, for example,~\cite{Chien2018,Ruggeri2023,Zhen2023}, among others). However, few of them present a topological or geometric perspective, favoring probabilistic or statistical approaches. We close the section reviewing the only two---to the best of our knowledge---geometric-based clustering method for hypergraphs: an extension of the Normalized Cut algorithm \cite{ShiMalik2000, MeilaShi2001, NgJordanWeiss2001} for spectral clustering in graphs \cite{ZhouHuangScholkopf2006,YangDengLu2021}; and an application of spectral clustering based on the Ollivier-Ricci curvature \cite{CoupetteDalleigerRieck2023}.

\subsection{Preliminaries on Hypergraphs}
\label{sec:hypergraph-cores-prelim}

Given a hypergraph $H=(V,E)$, we let  $n=|V|$, $m=|E|$, and $p=\sum_{e\in E}|e|$, where $|e|$ denotes the hyperedge size, that is, the number of vertices (or nodes) it contains.
We let the rank of the hypergraph $r$ to be the maximum hyperedge size. A hypergraph is \emph{capacitated} if it is equipped with a capacity function $c:E\to \mathbb{R}^+$, where $c(e)$ assigns a positive real number to each hyperedge $e \in E$, representing the capacity of that hyperedge. If all capacities are equal to one, the hypergraph is called \emph{uncapacitated}.
The \emph{$(k,q)$-core} of a hypergraph, defined by Ahmed et al. in~\cite{AhmedBatageljFu2007}, is a maximal subgraph in which each node has at least $k$ hypergraph degrees and each hyperedge contains at least $q$ nodes~\cite{LeeGohLee2023}.

Several notions of hypergraph Laplacians exist in the literature, each motivated by different applications or structural properties~\cite{AgarwalBransonBelongie2006,HeinSetzerJost2013,LiMilenkovic2018,SaitoMandicSuzuki2018,SaitoHerbster2023}. Here, we present one such formulation, based on the submodular Laplacian framework~\cite{Yoshida2019}, which is commonly used in recent works on spectral sparsification. 
The Laplacian $L_H \colon \mathbb{R}^V \rightarrow \mathbb{R}^V$ of an undirected weighted hypergraph $H$ with weight function $w \in \mathbb{R}^E_+$ is defined~\cite{Louis2015,Yoshida2019} such that its quadratic form $\mathbf{x}^{\top} L_H \mathbf{x}$ satisfies
\begin{equation*}
\mathbf{x}^{\top} L_H \mathbf{x} = \sum_{e \in E} 
w(e) \max _{u, v \in e}(\mathbf{x}(u)-\mathbf{x}(v))^2,
\end{equation*}
for every $\mathbf{x}\in \mathbb{R}^V$.

We say that a node set $S \subseteq V$ \emph{cuts} a hyperedge $e$ if $e\cap S \neq \emptyset$ and $e \cap (V \setminus S)\neq \emptyset$. Specifically, for a set $S \subseteq V$, the quadratic form $\mathbf{1}_S^{\top} L_H \mathbf{1}_S$ coincides with the cut size of $S$, where $\mathbf{1}_S \in \mathbb{R}^V$ is the characteristic vector of $S$.

A directed hypergraph $H = (V, E, w)$ as defined by Gallo et al. in~\cite{GalloLongoPallottino1993} consists of a node set $V$, a set of hyperedges $E$, and a weight function $w \in \mathbb{R}^E_+$, where each hyperedge $e$ is a pair $(T_e, H_e)$ of (not necessarily disjoint) sets of nodes, where $T_e,H_e \subseteq V$ are called the tail and head of $e$, respectively. 
The Laplacian $L_H \colon \mathbb{R}^V \rightarrow \mathbb{R}^V$ for a directed hypergraph $H$ can be derived using the submodular Laplacian framework. 
Using the submodular form of the Laplacian~\cite{Yoshida2019}, Soma and Yoshida~\cite{SomaYoshida2019} derived the quadratic form $\mathbf{x}^{\top} L_H \mathbf{x}$ for a directed hypergraph to be as follows:  
\[\mathbf{x}^{\top} L_H \mathbf{x} = \sum_{e = (T_e, H_e) \in E} w(e) \max_{u \in T_e} \max_{v \in H_e} ([\mathbf{x}(u) - \mathbf{x}(v)]^+ )^2,
\] 
where $[x]^+ = \max {( x, 0 )}$.

A weighted subgraph $H'$ of a hypergraph $H$ on a node set $V$ is said to be an $\varepsilon$-spectral sparsifier of $H$ if the following holds for every $x \in \mathbb{R}^V$~\cite{SomaYoshida2019}\footnote{Notice that there is a difference between the formulation from~\cite{SomaYoshida2019} vs. the formulation from~\cite{SpielmanSrivastava2011}.}, 
\[
(1 - \varepsilon)\mathbf{x}^{\top} L_{H'} \mathbf{x} \leq \mathbf{x}^{\top} L_H \mathbf{x} \leq (1 + \varepsilon)\mathbf{x}^{\top} L_{H'} \mathbf{x}. 
\]
\subsection{Percolation Methods for Finding Hypergraph Cores}
\label{sec:hypergraph-cores-percolation}

Percolation processes describe the size of the largest connected component when nodes or hyperedges are randomly removed~\cite{BianconiDor2024}.
In the context of studying the dependence of hyperedges on their nodes, also known as hypergraph percolation, 
as noted by Bianconi and Dor in~\cite{BianconiDor2024}, two types of hyperedges capture different classes of higher-order interactions. The first type, found in networks such as social interactions, remains intact unless all but one of its nodes are removed, or fail, meaning that the hyperedge collapses only when it loses all but one of its constituent nodes. This failure indicates that the functionality or connection of the hyperedge is lost when it is no longer supported by a sufficient number of nodes.
Thus, these hyperedges can withstand the failure of one or more nodes. The second type consists of hyperedges that fail once one of their nodes is damaged, such as in supply chains or protein-interaction networks. 
In the first category, Bianconi and Dorogovtsev~\cite{BianconiDor2024} noted that current theories and models that treat hypergraphs as factor graphs---which represent hypergraphs as bipartite networks of nodes and hyperedges---are completely effective.

We start by outlining contributions to the theory of percolation for hypergraphs with hyperedges of varying cardinalities. Sun and Bianconi~\cite{SunBianconi2021} studied \emph{hypergraph pruning} (cf., link pruning from \cref{sec:graph-core-FR-community}) at the node and hyperedge levels separately, proposing percolation methods for finding the $(k,1)$- and $(1,k)$-cores of random hypergraphs, respectively; this percolation procedure extends percolation on factor graphs. Higher-order percolation on hypergraphs is related to percolation on multiplex networks, that is, networks with multiple layers sharing the same set of nodes.   
Random hypergraphs can exhibit complex multiplex topologies, which are characterized by layers representing hyperedges of specific cardinalities, leading to the concept of random multiplex hypergraphs~\cite{SunBianconi2021}.  
This approach avoids the strict assumption of fixed cardinality hyperedges. The authors argued that multiplex hypergraphs are excellent tools in statistical mechanics for examining a wide range of higher-order percolation processes, applicable across various fields. They defined ensembles of random multiplex hypergraphs, where each node has a generalized degree vector representing the number of hyperedges of different cardinalities incident to it. These hypergraphs show significant interlayer hyperdegree correlations and differ from multiplex bipartite networks, although mapping between the two is possible under certain conditions. The study demonstrates that the complex topology and interlayer correlations of multiplex hypergraphs influence percolation processes, allowing adjustment of the percolation threshold. 
The research reveals that the multiplex nature of hypergraph ensembles can be used to explore higher-order percolation problems, uncovering complex behaviors such as discontinuous hybrid transitions and multiple percolation transitions.

Lee et al.~\cite{LeeGohLee2023} extended this concept by defining the $(k,q)$-core of a hypergraph, which is obtained by iteratively removing nodes with degrees less than $k$ and hyperedges smaller than $q$, leading to a repeated pruning process. 
As the process continues, any nodes and hyperedges that subsequently fall below these thresholds are also removed. This process results in a core subhypergraph where every node has at least $k$ connections, and every hyperedge connects to at least $q$ nodes.
To describe this pruning process for uncorrelated hypergraphs, that is, hypergraphs where there is no correlation between node degrees and hyperedge sizes, the authors derived evolution equations for the degree and size distributions. These equations, framed within a bipartite network representation of hypergraphs, track the fractions of nodes and hyperedges pruned at each step. The self-consistency equations for the core structure were obtained using generating functions, which facilitated the analysis of 
$(k,q)$-core percolation transitions. 

The analysis of Lee et al.~revealed that there is a hybrid phase transition of $(k,q)$-core percolation when $k\geq3$ or $q\geq 3$ and the second moments of the degree and size distributions are finite. This hybrid phase transition is characterized by a discontinuous jump in the order parameter, accompanied by critical behaviour at the transition point. Additionally, the relaxation dynamics of the fractions of nodes with degree 
$z=k-1$ and hyperedges of size $n=q-1$ at the transition point exhibit a universal behaviour proportional to $t^{-2}$ for $k\geq 3$ or $q \geq 3$, where $t$ denotes the time step.
When $k=q=2$, Lee et al.~analytically derived and numerically confirmed a novel degree-dependent critical relaxation dynamics, showing that $P^{(u)}(z,t) \sim t^{-z}$ and $P^{(e)}(n,t) \sim t^{-n} $ for $z$ and $n\geq 2$, with both decaying as $t^{-3}$ for $z=n=1$. Here, $P^{(u)}(z,t)$ represents the probability distribution of nodes with degree  $z$ at time $t$, and $P^{(e)}(n,t)$ represents the probability distribution of hyperedges with size $n$ at time $t$. 

Lee et al.~further demonstrated that 
$(k,q)$-core decomposition is more effective in eliminating redundant modular structures compared to traditional 
$k$-core decomposition. This obervation was highlighted through a case study on a coauthorship hypergraph, where 
$(k,q)$-core decomposition successfully identified high-impact teams by eliminating modular structures with large hyperedges and many low-degree nodes. Lee et al.~also noted that $(k,q)$-core decomposition has advantages in classifying groups (hyperedges) based on their properties more effectively than the previous 
$k$-core decomposition via this case study. 
Whereas the theoretical analysis was confined to hypergraphs with homogeneous distributions, the authors remarked that their theory can be generalized to heterogeneous cases due to the convergence of the infinite series, irrespective of the moments' divergence in the degree distribution. Furthermore, the method allows for extensions to cases with degree-size correlation. A key shortcoming of the existing theory, as pointed out in~\cite{LeeGohLee2023}, is its reliance on a tree-like assumption. 
This relaxation and the development of broader theoretical models represent promising directions for future research.

Percolation theory for hypergraphs with hyperedges belonging to the second type, i.e.,~hyperedges that fail once one of their nodes is damaged, was recently developed by Bianconi and Dorogovtsev~\cite{BianconiDor2024}. Whereas hypergraphs can be efficiently represented by factor graphs, $k$-core percolation can significantly differ between the two. 
To address these disparities, Bianconi and Dorogovtsev introduced a set of pruning processes targeting either nodes or hyperedges based on connectivity in the second neighbourhood. 
They constructed a message-passing theory for hypergraph percolation, utilizing the generating function formalism and supporting their analysis with Monte Carlo simulations on both random and real-world hypergraphs. They further derived message-passing equations for percolation on both factor graphs and hypergraphs. Next, they applied message-passing theory  to random hypergraphs described by a degree distribution for nodes and a cardinality distribution for hyperedges, and a joint degree distribution where each degree is a vector listing the numbers of hyperedges of each cardinality adjacent to a node. 

Bianconi and Dorogovtsev further investigated  the critical behavior of hypergraph percolation, revealing that the node percolation threshold for hypergraphs is higher than that for factor graphs. Additionally, unlike ordinary graphs where node and edge percolation thresholds coincide, in hypergraphs the situation is different, with the threshold for node percolation being higher than that for hyperedge percolation. The authors determined the presence of a percolation cluster (i.e.,~a giant connected component) and its relative size through self-consistent equations, which are derived using generating functions. The study shows that in hypergraphs, distributions of hyperedge sizes with heavy tails do not result in hyperresilience, which is a contrast to factor graphs. In factor graphs, if the size distribution has a diverging second moment, it leads to a percolation threshold of zero. The distinction between node and hyperedge percolation in hypergraphs is significant, especially when hyperedges have large cardinalities. Notably, when pruning focuses solely on hyperedges, the phase diagram converges to that of factor graph $k$-cores. 

Importantly, the message-passing algorithm developed by Bianconi and Dorogovtsev does not assume the absence of correlations in hypergraphs, which allows the problem to be treated numerically. The final formulas are obtained for hypergraphs without degree-degree correlations between different nodes. The authors suggested that the case of correlated hypergraphs can also be treated analytically within their framework.

Recently, Bianconi and Dorogovtsev~\cite{BianconiDoro2024} studied $k$-core percolation processes on hypergraphs and contrasted these with the analogous processes on factor graphs. The critical distinction arises in the integrity of hyperedges: in hypergraphs, every node within a hyperedge must be intact for the hyperedge to be considered intact, reflecting the all-or-none characteristic observed in systems such as supply chains, protein-interaction networks, and chemical reactions. The authors formulated a message-passing theory tailored to the $k$-core percolation on hypergraphs, leveraging the theory of critical phenomena on networks to underline the main differences from factor graph-based percolation where hyperedges can still function despite partial node failures. They explored phase transitions and phase diagrams in these systems, noting more complex behaviors in hypergraphs compared to factor graphs, especially in scenarios where node or hyperedge pruning is based on broader neighbourhood connectivity. To resolve discrepancies between $k$-core percolation on hypergraphs and factor graphs, they introduced a series of pruning processes, designed to either exclusively target nodes or hyperedges, with a dependence on their connectivity to the second neighbourhood. 

\subsection{Spectral Sparsification Methods}
\label{sec:hypergraph-cores-spectral}

Soma and Yoshida~\cite{SomaYoshida2019} extended spectral sparsification of graphs to hypergraphs by providing an $\varepsilon$-spectral sparsifier for a hypergraph $H$.  
They presented a polynomial-time algorithm that constructs an $\varepsilon$-spectral sparsifier of an undirected or directed hypergraph $H=(V,E)$ with $O(n^3 \log n / \varepsilon^2)$ hyperedges, where $n=|V|$.
Their main contributions include the development of a randomized algorithm for spectral sparsification and its application to various computational problems involving hypergraph Laplacians. For an undirected hypergraph, their algorithm outputs an $\varepsilon$-spectral sparsifier $H$ with $O(n^3 \log n / \varepsilon^2)$ hyperedges with high probability in $O(pn + m \log(1/\varepsilon^2) + n^3 \log n / \varepsilon^2)$ time, where $n=|V|$, $m=|E|$, $p=\sum_{e \in E}|e|$.  

Soma and Yoshida further showed that spectral sparsification of directed hypergraphs can also be achieved with $O(n^3 \log n / \varepsilon^2)$ hyperedges. This result is significant because even directed graphs with $O(n^2)$ edges do not admit nontrivial cut/spectral sparsification of size $o(n^2)$ ($o(\cdot)$ representing little-$O$ notation), whereas directed hypergraphs, which could have $O(4^n)$ hyperedges, do admit such sparsification. Their study also demonstrated how spectral sparsification can be applied to enhance the time and space  efficiency of algorithms that work with the quadratic forms of hypergraph Laplacians. This study includes tasks such as calculating eigenvalues, solving systems based on Laplacians, and conducting semi-supervised learning. 

Kapralov et al.~\cite{KapralovKrauthgamerTardos2021} presented a polynomial-time algorithm that constructs an $\varepsilon$-spectral sparsifier for a hypergraph with only $O^*(nr)$ hyperedges, where $r$ is the maximum hyperedge size, and the notation $O^*(\cdot)$ suppresses $(\varepsilon^{-1}\log n)^{O(1)}$ factors~\cite{KapralovKrauthgamerTardos2021}. This result improves upon the previous bounds, i.e.,~$O^*(n^3)$ set by Soma and Yoshida~\cite{SomaYoshida2019} 
and $O^*(n r^3)$ by Bansal et al.~\cite{BansalSvenssonTrevisan2019}. 
The algorithm achieves this sparsification with high probability and a running time of $O(m r^2) + n^{O(1)}$. 
Additionally, the authors established lower bounds on the bit complexity for any compression scheme that approximates all cuts in a hypergraph within a factor of $1 \pm \varepsilon$. We refer the reader to \cref{sec:hypergraph-cores-prelim} for the definition of a hypergraph cut. They demonstrated that $\Omega(n r)$ bits are required to represent all cut values of an $r$-uniform hypergraph where $r$ equals $n^{O(1 / \log \log n)}$. 
For directed hypergraphs, Kapralov et al. further presented an algorithm that computes an $\varepsilon$-spectral sparsifier with $O^*(n^2 r^3)$ hyperedges, which represents a significant improvement over the previously known $O^*(n^3)$ bound by Soma and Yoshida. This result is particularly impactful when the rank $r$ is small. They also introduced a new approach for proving concentration of the nonlinear version of the quadratic form associated with the Laplacians in hypergraph expanders.

Independently, Lee \cite{Lee2023} and Jambulapati et al.~\cite{JambulapatiLiuSidford2023} further improved the bound in the undirected hypergraph with a near-linear time algorithm. They introduced an algorithm that constructs an $\varepsilon$-spectral sparsifier with $O(n \varepsilon^{-2} \log n \log r)$ hyperedges for an $n$-node, $m$-edge hypergraph of rank $r$. This algorithm operates in nearly linear $O^*(mr)$ time, significantly enhancing both the size and efficiency over previous works (e.g.,~\cite{BansalSvenssonTrevisan2019,KapralovKrauthgamerTardos2021}). 
The algorithm by Jambulapati et al.~\cite{JambulapatiLiuSidford2023} utilizes a novel approach involving \emph{group leverage score} overestimates and generic chaining, enhancing the calculation of sampling weights in nearly linear time to achieve optimal sparsification bounds. In the context of hypergraphs, group leverage scores assess the importance or influence of hyperedges, typically derived from spectral properties. By overestimating these scores, the algorithm simplifies the calculations, using approximations to determine the importance of each hyperedge efficiently. This approximation results in adjusted sampling probabilities for retaining significant hyperedges in the sparsifier, which maintains the structural integrity of the original hypergraph.

Furthermore, the technique of generic chaining, developed by Talagrand~\cite{Talagrand2014}, was adapted to manage the accumulation of error in sampling. This method constructs multiple layers of increasingly finer nets over the probability space, chaining these layers to control the maximum deviation of the stochastic processes involved. This rigorous control allows for a more refined process in determining which hyperedges to retain, tightening the bound on the number of necessary hyperedges. These innovations offer  improvements over previous methods, demonstrating that every hypergraph admits a sparsifier with nearly linear time complexity and fewer hyperedges. Specifically, this method effectively reduces the bound on the number of hyperedges to $O(n \varepsilon^{-2} \log n \log r)$, improving over the previous bounds and approaching the theoretical limits of sparsification.

Given a hypergraph $H=(V,E)$ with $n=|V|$, $m=|E|$, and $p=\sum_{e \in E}|e|$, Chekuri and Xu~\cite{ChekuriXu2018} extended sparsification techniques from graphs to uncapacitated hypergraphs, focusing on constructing a $k$-trimmed certificate.
 A \emph{$k$-trimmed certificate} is a (sparse) subhypergraph $H'$ obtained via hyperedge deletion and trimming that preserves all cuts of value up to $k$~\cite{ChekuriXu2018}. 
The authors improved on previously bounds in~\cite{queyranne1998minimizing,mak2000fast,klimmek1996simple} to ﬁnd the min-cut (or minimum cut) of an
arbitrary symmetric submodular function.  
For a hypergraph, their method adapts the \emph{maximum adjacency ordering} to arrange nodes in a sequence where each node has high connectivity to previously ordered nodes. Edges are then evaluated, and for each node $v$, only the first $k$ significant backward edges (i.e.,~those connecting 
$v$ to previously ordered nodes) are retained. 
This process trims the hypergraph by removing less critical connections while maintaining essential connectivity properties. The construction of the $k$-trimmed certificate involves creating a data structure in $O(p)$ time. The practical implications include significantly faster algorithms for computing the global min-cut in hypergraphs. For example, the global min-cut can be computed in $O(p+\lambda n^2)$ time, where $\lambda$ is the min-cut value of the hypergraph. 
Chekuri and Xu proposed a split oracle that, given a hypergraph $H$ and its min-cut value $\lambda$, determines if a split exists between any two nodes $\{s,t\}$. The oracle operates in near-linear time by converting $H$ into an equivalent directed graph $H_D$ and computing the maximum $s$-$t$ flow in $H_D$. If this flow exceeds $\lambda$, no $s$-$t$ split exists; otherwise, a nontrivial min-$s$-$t$ cut implies a split. Canonical and prime decompositions are introduced to further break down a hypergraph. A hypergraph is \emph{prime} if it contains no splits, with all min-cuts being trivial. A \emph{canonical} decomposition captures all min-cut information and is constructed by refining the hypergraph iteratively using splits. This decomposition helps in identifying cores, as they correspond to components maintaining the hypergraph connectivity. To compute the canonical decomposition, an algorithm starts with a prime decomposition and iteratively glues components. This algorithm operates in $O(np + n^2 \log n)$ time for capacitated hypergraphs and $O(np)$ time for uncapacitated hypergraphs. 

Chekuri and Xu~\cite{ChekuriXu2018} discussed a compact representation of a hypergraph called the \emph{hypercactus representation}, which captures all min-cuts and is derived from the canonical decomposition. For a hypergraph with $\lambda = 1$, the hypercactus is constructed by identifying marker nodes and combining canonical decomposition components. The hypercactus representation can be computed in $O(n(p + n \log n))$ time and $O(p)$ space for capacitated hypergraphs, and in $O(p + \lambda n^2)$ time for uncapacitated hypergraphs. The authors then extended Matula's algorithm to hypergraphs, obtaining a $(2+\varepsilon)$-approximation for the global min-cut of a capacitated hypergraph in $O(\frac{1}{\varepsilon}(p \log n + n \log^2 n))$ time and for uncapacitated hypergraphs in $O(p/\varepsilon)$ time. This generalized version of Matula's algorithm incorporates the use of node orderings, inspired by Nagamochi and Ibaraki's~\cite{NagamochiIbaraki1992} maximum adjacency ordering method for graphs. Unlike graphs, hypergraphs have several possible orderings, and these yield different insights into their structure. The authors present an algorithm that computes approximate strengths for all edges in a hypergraph, running in $O(p \log^2 n \log p)$ time. The \emph{strength} $\gamma_H(e)$ of an edge $e$ is defined as the maximum min-cut value over all node-induced subhypergraphs containing $e$. The algorithm ensures that the estimated strength $\gamma'(e)$ is less than or equal to the actual strength $\gamma_H(e)$, and that the sum of the inverses of the approximate strengths satisfies certain bounds. This edge strength estimation facilitates the construction of a $(1+\varepsilon)$-cut sparsifier, a reduced hypergraph that preserves cut size within a factor of $1 \pm \varepsilon$. This cut sparsifier can be found in near-linear time, significantly reducing the number of edges while maintaining essential connectivity properties. 
This approach extends the work of Kogan and Krauthgamer~\cite{KoganKrauthgamer2015} and leads to faster algorithms for solving various cut and flow problems in hypergraphs of small rank. By using the approximate strengths to sample edges and create sparsifiers, the authors achieved  efficient cut approximations.

Chen et al. \cite{ChenKhannaNagda2020} improved on the work of Chekuri and Xu, providing an algorithm that runs in $\tilde{O}(mn + n^{10}/\varepsilon^7)$ time that constructs an approximate sparsifier of size $O\left(\varepsilon^{-2}n\log n\right)$. The authors claimed that this size bound is the best possible bound within a logarithmic factor by providing an example of a hypergraph such that any sparsifier must include every hyperedge. Their algorithm works by constructing a sparsifier via an auxiliary graph $G$ as follows. First, for each hyperedge $e$ in $H$, add a clique $F_e$ to $G$ whose node set is the same as the nodes of $e$, then assign weights to the edges in $G$ according to a balanced weight scheme, where the weights of the edges in $F_e$ range from $0$ to the weight of the hyperedge $e$. With this setup, the probability of sampling a hyperedge $e$ in $H$ is determined by the balanced weight scheme on the clique $F_e \in G$.

By generalizing the notion of effective resistance to hypergraphs, Aghdaei and Feng \cite{AghdaeiFeng2022} introduced a scalable algorithmic framework, called HyperEF, for spectral coarsening of large-scale hypergraphs.
This generalization relies on the nonlinear quadratic form introduced by Chan et al. in \cite{ChanLouieTang2018}: for a weighted undirected hypergraph $H=(V,E,w)$ with non-negative weights, and a vector $\mathbf{x}\in \mathbb{R}^{|V|}$, 
\[
Q(\mathbf{x}):=w_e\sum_{e\in E}\max_{v_{k},v_{l}\in e}(\mathbf{x}(k)-\mathbf{x}(l))^2,
\]
where $\mathbf{x}(k)$ is the $k$-th element of the vector $\mathbf{x}$.
In the graph setting, the effective resistance between two nodes $v_i$ and $v_j$ satisfies the following property, as shown in Theorem 1 in \cite{AghdaeiFeng2022}:
\[
r_{v_iv_j}=\max_{\mathbf{x}\in \mathbb{R}^{|V|}}\frac{(\mathbf{x}^{\top}(\mathbf{e}_i-\mathbf{e}_j))^2}{\mathbf{x}^{\top}L\mathbf{x}}.
\]

By extending the above formula to the hypergraph setting while replacing $\mathbf{x}^{\top}L\mathbf{x}$ with the nonlinear quadratic form $Q(\mathbf{x})$, Aghdaei and Feng~\cite{AghdaeiFeng2022} defined the effective resistance of a hyperedge $e$ as 
\begin{equation}
\label{eq:ER-hyper}
    R_e :=\max_{\mathbf{x}\in \mathbb{R}^{|V|}}\frac{\left(\mathbf{x}^{\top}\left(\mathbf{e}_i-\mathbf{e}_j\right)\right)^2}{Q(\mathbf{x})},
\end{equation}
assuming $v_i$ and $v_j$ are some nodes within $e$.
To enhance computational efficiency, $R_e$ is computed approximately as follows.
For a vector $\mathbf{x} \in \mathbb{R}^{|V|}$, we call $\tilde{R}_e(\mathbf{x}):=\frac{(\mathbf{x}^{\top}(\mathbf{e}_i-\mathbf{e}_j))^2}{Q(\mathbf{x})}$ the resistance ratio associated with $\mathbf{x}$.
Then, $R_e$ is approximated as the sum of the resistance ratios associated with a few Laplacian eigenvectors $\{\xi_i\}$ of a certain bipartite graph $G_b$ converted from the original hypergraph. This approximation is motivated by the following property of graph effective resistance: assuming $\{\xi_i\}_{i=1}^{|V|}$ is an orthonormal eigenbasis of the graph Laplacian associated with eigenvalues $\{\lambda_i\}_{i=1}^{|V|}$, then
\[
r_{v_iv_j}=\sum_{\lambda_i\neq 0}\frac{(\mathbf{x}^{\top}(\mathbf{e}_i-\mathbf{e}_j))^2}{\xi_i^{\top}L\xi_i}.
\]
Here, the vectors $\{\xi_i\}$ are computed approximately as the orthonormal basis of the Krylov subspace $K:=\mathrm{span}\{x, Ax, \dots, Ax^{\rho}\}$, for a non-negative integer $\rho$, a random vector $x$ and the adjacency matrix $A$ of the bipartite graph $G_b$.
Nodes $v_i$ and $v_j$ in \cref{eq:ER-hyper} are selected as the maximally separated nodes within this $\rho$-dimensional embedding space $K$.

The HyperEF algorithm iteratively contracts the hyperedges and nodes via the following steps: (1) compute the approximate effective resistance as described above; (2) compute the node weights via a specified formula; (3) compute the hyperedge weights using both the approximate effective resistance and the node weights; and (4) contract the hyperedges and cluster nodes when their weights are below a threshold. 
The performance of the algorithm is tested by its preservation of the so-called average conductance of clusters.
Experiment results on real-world VLSI designs (i.e., very-large-scale integration, the process of creating an integrated circuit by combining thousands to millions of transistors on a single chip) showed that HyperEF is faster than HMetis~\cite{karypis1998hypergraph}, a well-known hypergraph partitioner, while maintaining comparable or better average conductance values in most cases.

\subsection{Spectral Clustering for Hypergraphs}

A successful clustering technique in the graph domain is \emph{spectral clustering}, based on the eigenvalues and eigenvectors of some matrix containing geometric information about the graph. A well-established method for spectral clustering is the normalized cut algorithm \cite{ShiMalik2000,MeilaShi2001}, which intuitively minimizes the weights of the edges connecting different clusters without favoring clusters of isolated nodes in the graph. 

This normalized cut algorithm was extended to hypergraphs by Zhou et al.~\cite{ZhouHuangScholkopf2006}. The authors showed that, much like in the graph case, a relaxation in the hypergraph cut criterion allows to express the optimal cut from a decomposition of a positive semidefinite matrix, which can be regarded as a hypergraph Laplacian. Nonetheless, their approach had two limitations. On the one hand, it assumed that nodes within the same hyperedge have equal relevance, which can be restrictive. On the other hand, the need to solve an eigenproblem made the approach time and storage consuming.

GraphLSHC \cite{YangDengLu2021} was proposed by Yang et al.~as an extension of the previous method to overcome such limitations. For the importance of edges problem, the authors proposed to use weighted hypergraphs, and extended the normalized cut algorithm to such a framework. Here, the weights refer to ``intra-weights'' within hyperedges.  To address the computational issues, they proposed several optimizations, including an ``eigen-trick'' that is particularly designed to accelerate the eigenproblem solver.

To the best of our knowledge, the only existing curvature-based spectral clustering algorithm for hypergraphs was proposed by Coupette et al.~\cite{CoupetteDalleigerRieck2023}, who proposed an flexible framework to extend the Ollivier-Ricci curvature to hypergraphs (see \cref{sec:OR}). As part of their experimental exploration on \emph{hypergraph learning}, they performed spectral clustering using curvatures and other local edge features (e.g.,~the number of neighbours of a given edge) and node features (e.g.,~the average curvature of the incident edges in a node, or the size of its neighbourhood). They evaluated the quality of the clustering through the \emph{Wasserstein clustering coefficient}, which measures average intra-cluster Wasserstein distances to inter-cluster Wasserstein distances; and observed that curvature-based clustering consistently obtains better results.  

\begin{table}[!t]
\caption{A table summarizing recent works in finding cores of simplicial complexes, as surveyed in \cref{sec:sc-cores}. 
Ref.: references. SC: simplicial complex. WSSD: well-separated simplicial decomposition. PD: persistence diagram. VR: Vietoris--Rips complex. \v{C}ech: \v{C}{e}ch complex. ER: effective resistance.}
\label{table:sc-core-summary}
\centering
\begin{tabular}{p{1cm}p{1cm}p{4cm}p{2.7cm}p{2.5cm}}
\hline\noalign{\smallskip}
\textbf{Ref.} & \textbf{SC} & \textbf{Key component(s)} & \textbf{Criteria} & \textbf{Property Preserved} \\
\noalign{\smallskip}\svhline\noalign{\smallskip}
\cite{Sheehy2012}            & VR         & point collapse/removal                  & metric ball cover   & PD \\
\cite{DeyFanWang2014}        & VR         & batch collapse with pairwise distance   & metric ball cover   & PD \\
\cite{DeyShiWang2019}        & VR         & batch collapse with set distance        & metric ball cover   & PD \\
\cite{KerberSharathkumar2013}& \v{C}ech   & WSSD                                    & metric ball cover   & PD \\
\cite{MemoliOkutan2021}      & VR         & vertex quasi-distance                   & codensity           & PD \\
\cite{BauerEdelsbrunner2017} & SC         & simplicial collapse                     & none                & PD \\
\noalign{\smallskip}\hline\noalign{\smallskip}
\cite{CaoMonod2022}          & VR         & subsample smaller point clouds          & none                & PD \\
\cite{GomezMemoli2024}       & VR         & subsample smaller point clouds          & none                & PD \\
\cite{SolomonWagnerBendich2022} & VR      & subsample smaller point clouds          & none                & PD \\
\noalign{\smallskip}\hline\noalign{\smallskip}
\cite{OstingPalandeWang2020} & SC         & Sampling with replacement               & ER                  & spectral \\
\noalign{\smallskip}\hline\noalign{\smallskip}
\cite{EbliSpreemann2019}     & SC         & harmonic embedding, subspace clustering & none                & homology generators (experimentally) \\
\noalign{\smallskip}\hline
\end{tabular}
\end{table}

\section{Finding Cores of Simplicial Complexes}
\label{sec:sc-cores}

Finally, we review recent works in finding cores of simplicial complexes. Majority of these works arise from the fields of computational topology; thus, they have a strong focus on preserving homological information of the complexes. 
After reviewing various notions of simplicial complexes in \cref{sec:sc-cores-prelim},  
we discuss sparsification methods that preserve (persistent) homology (\cref{sec:sc-cores-PH,sec:sc-cores-sampling}) and spectral properties (\cref{sec:sc-cores-spectral}) of the complexes.  
We finish this section with a discussion on simplicial spectral clustering (\cref{sec:sc-cores-clustering}).
Whereas there have been recent works on percolation theory for simplicial complexes, e.g., \emph{homological percolation} by Bobrowski and Skraba~\cite{BobrowskiSkraba2020}, different from hypergraphs, these works have yet to produce well-defined cores. Therefore, in this paper, we choose not to survey simplicial percolation. 

\subsection{Preliminaries on Simplicial Complexes}
\label{sec:sc-cores-prelim}

We first review the definitions of various simplicial complexes popular in computational topological, including \v{C}ech, Vietoris--Rips, Delaunay, Alpha, Delaunay-\v{C}ech, and Wrap complexes.  

Given a finite set $X \subseteq \Rspace^n$ and a radius $r \geq 0$,  
the \emph{\v{C}ech complex} of $X$ is the nerve of the $r$-balls centered at points of $X$, that is,  
\[
\Cech_r(X) = \{Q \subseteq X \mid \bigcap_{x \in Q} B_x(r) \neq \emptyset\}. 
\]
The Vietoris--Rips complex for a scale $r$ is 
\[
\VR_r(X) = \{Q \subseteq X \mid \diam(Q) \leq r\}.
\]
The \emph{Voronoi cell} of a point $x \in X$ (w.r.t. $X$) contains all points closer to $x$ than to any other point, 
\[
\Vor(x) = \{y \in \Rspace^n \mid d(y,x) \leq d(y,p), \; \forall \; p \in X\}. 
\]
The \emph{Voronoi ball} of $x$ (w.r.t. $X$) for a radius $r$ is 
\[
\Vor_r(x) = B_r(x) \cap \Vor(x). 
\]
The \emph{Delaunay complex} of $X$ for a radius $r$ is 
\[
\Del_r(X) = \{Q \subseteq X \mid \bigcap_{x \in Q} \Vor_r(x) \neq \emptyset \}.
\]
$\Del_r(X)$ is also called the \emph{alpha complex}.
The \emph{Delaunay triangulation} of $X$ is $Del(X):=Del_\infty(X)$.

The \emph{Delaunay-\v{C}ech} complex for a radius $r$ restricts the \v{C}ech complex to the Delaunay triangulation~\cite{BauerEdelsbrunner2017}:
\[
\DelCech_r(X)=\{Q \in \Del(X) \mid \bigcap_{x \in Q} B_r(x) \neq \emptyset\}. 
\]
The \emph{Wrap complex} is defined using the gradient of the Delaunay radius function $s_X: \Del(X) \to \Rspace$, see~\cite{BauerEdelsbrunner2017} for its technical definition.   

\subsection{Persistent Homology}
\label{sec:persistence}

In this section, we introduce the basics pertaining to the theory of persistence that are necessary in~\cref{sec:sc-cores}. 
The fundamental idea is to build a geometric representation of the  input data in the form of a \emph{filtration}, from which one computes invariants capturing the topological features at multiple scales. 
For the purpose of this survey, given a simplicial complex $K$, a filtration is an increasing sequence of simplicial complexes connected by inclusions, 
\[
\emptyset = K_0 \hookrightarrow K_1 \hookrightarrow \dots \hookrightarrow K_N = K.  
\]
To extract the topological information from a filtration, persistent homology studies the \(p\)-dimensional homology groups $\Hgroup_p$ of the simplicial complexes in the filtration. 
For our purposes, these are vector spaces that contain information about the topological features in a simplicial complex: connected  components (for \(p=0\)), tunnels (for \(p=1\)), voids (for \(p=2\)), and so on. 
By computing $p$-homology groups for all simplicial complexes in a filtration, one obtains a family of vector spaces connected by linear maps that satisfy a commutativity relation, 
\[
\emptyset = \Hgroup_p(K_0) \to \Hgroup_p(K_1) \to \dots \to \Hgroup_p(K_N) = \Hgroup_p(K).  
\]
This structure is an example of a \emph{persistence module}. 
Algebraically, a persistence module, denoted as $(M, \phi)$, is a functor $M: \mathbf{P} \to \mathbf{Vec}$, where $\mathbf{P}$ is a poset category (e.g., $\mathbf{P} = \Rspace^n$), and $\mathbf{Vec}$ is the category of finite-dimensional vector space over a field~\cite{Bjerkevik2021}. The \emph{transition map} between a pair of such vector spaces is a linear transformation $\phi_{s}^{t}: M_s \to M_t$ for $s \leq t \in \mathbf{P}$.   
When $\mathbf{P}=\Rspace$, by the Structure Theorem \cite{ZomorodianCarlsson2005, Crawley2015, BotnanCrawley2020}, a persistence module is fully explained by a family of intervals called the \emph{persistence barcode} of the module. 
An interval or \emph{bar} in a barcode, denoted as $(x,y)$, represents the birth time $x$ and the death time $y$ of a homological feature in the filtration. 
An equivalent representation of the barcode is given by the \emph{persistence diagram}, a multiset of points in the extended real plane $(\Rspace \cup {\pm\infty})^2$.
A point $(x,y)$ in the persistence diagram corresponds to a bar in the barcode. 
The persistence diagram has been extensively investigated and successfully applied to study sensor networks~\cite{deSilvaGhrist2007},
protein interactions~\cite{GameiroHiraokaIzumi2014,KovacevNikolicBubenikNikolic2016,XiaLiMu2018}, DNA structures~\cite{EmmettSchweinhartRabadan2016}, robot trajectories~\cite{PokornyHawaslyRamamoorthy2016}, and bipedal walks~\cite{VasudevanAmesBajcsy2013}, to name a few.  

The interleaving distance describes the proximity between two persistence modules. 
Two persistence modules, $(M,\phi)$ and $(N, \psi)$, are \emph{$\varepsilon$-interleaved} if there are families of functions \(f_t: M_t \to N_{t+\varepsilon}\) and \(g_t: N_t \to M_{t+\varepsilon}\) that commute with the internal transition maps, i.e., the following diagrams commute for any $s, t \in \Rspace$ where $s \leq t$:   
\[
\begin{tikzcd}
M_s \arrow[rr, "{\phi_{s,t}}"] \arrow[rd, "f_s"] &                                                                  & M_t \arrow[rd, "f_t"] &                \\
                                                      & N_{s+\varepsilon} \arrow[rr, "{\psi_{s+\varepsilon,t+\varepsilon}}"] &                       & N_{t+\varepsilon},
\end{tikzcd} 
\begin{tikzcd}
N_s \arrow[rr, "{\psi_{s,t}}"] \arrow[rd, "g_s"] &                                                                  & N_t \arrow[rd, "g_t"] &                \\
  & M_{s+\varepsilon} \arrow[rr, "{\phi_{s+\varepsilon,t+\varepsilon}}"] &                       & M_{t+\varepsilon}, 
\end{tikzcd}\]
\[
\begin{tikzcd}
M_t \arrow[rd, "f_t"] \arrow[rr, "{\phi_{t,\, t+2\varepsilon}}"] &                                             & M_{t+2\varepsilon}, \\
 & N_{t+\varepsilon} \arrow[ru, "g_{t+\varepsilon}"] &                
\end{tikzcd} \quad
\begin{tikzcd}
N_t \arrow[rd, "g_t"] \arrow[rr, "{\psi_{t,\, t+2\varepsilon}}"] &                                             & N_{t+2\varepsilon}. \\
 & M_{t+\varepsilon} \arrow[ru, "f_{t+\varepsilon}"] &                
\end{tikzcd}\]
The interleaving distance between two persistence modules is 
\[
\inf\{\varepsilon \mid (M,\phi) \text{ and } (N, \psi) \text{ are \emph{$\varepsilon$-interleaved}}\}.
\]

Given two filtrations $\Kcal:=(K_t)$, $\Lcal:=(L_{t})$, the persistence diagram of $\Kcal$ is a $c$-approximation~\cite{Sheehy2012} to that of $\Lcal$ if there is a bijection $\pi$ between the diagrams such that the birth (resp. death) time of a point $x$ in the persistence diagram of $\Kcal$ and $\pi(x)$ in the persistence diagram of $\Lcal$ differ by at most a factor of $c$.

\subsection{Sparse Filtrations Approximating Persistent Homology}
\label{sec:sc-cores-PH}

A significant issue in computational topology is that the number of simplices in a filtration may explode as the number of data points increases. 
For example, for large enough radii, both the \v{C}ech and the Vietoris--Rips complexes have $k$-skeletons (i.e., simplices up to dimension $k$) of size $O(n^{k+1})$ \cite{Sheehy2012}, and both complexes contain up to $2^n - 1$ (i.e.,~$O(2^n)$) simplices for $n$ data points (see Table 1 in \cite{otter2017roadmap}). 
Therefore, many works focus on sparsifying the simplicial filtrations  (i.e.,~removing points or simplices) while approximately preserving their homological information. We review  representative works based on obtaining sparse filtrations~\cite{Sheehy2012,CavannaJahanseirSheehy2015} and/or  applying batch collapse~\cite{DeyFanWang2014}. 
 
Given $n$ number of points, Sheehy~\cite{Sheehy2012} constructed a linear-size (i.e., $O(n)$) filtered simplicial complex whose persistence diagram approximates that of a Vietoris--Rips filtration.
The key idea is the introduction of a \emph{relaxed Vietoris--Rips filtration} based on the notion of a \emph{relaxed distance} that is provably close to the input metric. This new relaxed distance adds a weight to each point, which in turn shrinks the metric ball centered at the point. If a ball is covered by nearby balls, its center point can be deleted without changing the topology of the input point set, resulting in a sparsified complex. 
Formally, given a finite subset $P$ of a metric space $(X,d)$ and a user-defined parameter $\varepsilon \leq \frac{1}{3}$, the algorithm assigns a deletion time $t_x$ (based on a hierarchical net-tree construction, see~\cite{Sheehy2012} for details) to each point $x \in P$ and defines its weight $w_x$  based on a scale parameter $r$:  
\[w_x(r):= \begin{cases}
    0 & \text{if }r\leq (1-2\varepsilon)t_x;\\
    \frac{1}{2}(r - (1-2\varepsilon) t_x) & \text{if }(1-2\varepsilon)t_x < r< t_x; \\
    \varepsilon r & \text{if }t_x \leq r.
\end{cases}\]
Using these weights, the \emph{relaxed distance} at scale $r$ between two points $x,y \in P$ is 
\[ 
\hat{d}_r(x,y) := d(x,y) + w_x(r) + w_y(r).
\]
The \emph{relaxed Vietoris--Rips complex} is defined to be $\RVR_r:=\VR_r(P, \hat{d}_r)$, which gives rise to a filtration with varying $r$.  
Let $\Rcal_r := \VR_r(P,d)$. The multiplicative interleaving is 
\[\Rcal_{\frac{r}{c}}\subseteq \hat{\Rcal}_r \subseteq \Rcal_r,\] where $c = \frac{1}{1-2\varepsilon}$ (see Lemma 4.2.~in \cite{Sheehy2012}).  

The algorithm then constructs a \emph{sparse zigzag Vietoris--Rips filtration} with $O(n)$ total number of simplices, whose persistence diagram is identical to that of the relaxed Vietoris–Rips filtration.  
To obtain such a sparsified filtration, the algorithm deletes vertices and all the corresponding incident simplices from the relaxed Vietoris–Rips filtration, when they become unimportant at the scale being considered. 
The algorithm makes use of an \emph{open net} $\Ncal_r$ and a \emph{closed net} $\overline{\Ncal}_r$ at scale $r$, defined to be the points in $P$ with deletion time above $r$,   
\[
\Ncal_r := \{x \in P : t_x > r\};\;\;  
\overline{\Ncal}_r := \{x \in P : t_x \geq r\}. 
\]
The \emph{sparse zigzag Vietoris--Rips complex} at scale $r$ is the subcomplex of $\hat{\Rcal}_r$ induced by the vertices of $\Ncal_r$,
\[
\Qcal_r := \{\sigma \in \hat{\Rcal}_r: \sigma \subseteq \Ncal_r\}= \VR_r(\Ncal_r, \hat{d}_r).
\] 
Its closed version is defined similarly, 
$\overline{\Qcal}_r := \VR_r(\overline{\Ncal}_r, \hat{d}_r).$ 

The \emph{sparse zigzag Vietoris--Rips filtration} is 
\[
\dots \hookrightarrow \overline{\Qcal}_{r'} \hookleftarrow {\Qcal}_{r'} \hookrightarrow \overline{\Qcal}_{r} \hookleftarrow \dots 
\]
This leads to a filtration without the zigzag, defined by the union 
$\Scal_r = \bigcup_{r' \leq r} \overline{\Qcal}_{r'},$
which gives rise to the so-called \emph{sparse Vietoris--Rips filtration}, having the same persistence diagram as the zigzag one.

Sheehy provided theoretical guarantees that (i) the persistence diagram of a sparse Vietoris--Rips filtration is a multiplicative $c$-approximation (see \cref{sec:persistence} for the definition) of the standard one and (ii) the filtration has linear size~\cite{Sheehy2012}.   

Dey et al.~\cite{DeyFanWang2014} provided an efficient algorithm to compute persistence diagrams of filtrations where the simplicial complexes are connected by simplicial maps instead of inclusions (cf.,~\cref{sec:persistence}).  
They applied their algorithm to approximate the persistence diagram of a Vietoris--Rips filtration via a sparsified filtration. 
The key idea is that a simplicial map can be decomposed into  elementary inclusions and vertex collapses, and these atomic operations can be used to collapse input points in batches with an increasing radius, thereby controlling the size of the filtration during persistence computation. 

Using a mapping cylinder construction, a simplicial map $f: K_s \to K_t$ can be simulated with a zigzag connected by inclusions $K_s \hookrightarrow \hat{K} \hookleftarrow K_t$, via a third, significantly larger simplicial complex $\hat{K}$. 
Dey et al.~provided an improved construction that converts a zigzag module connected by simplicial maps into a zigzag module connected only by inclusions, where the intermediate modules are not as big as the naive construction of $\hat{K}$. 
They also showed that for a monotone sequence of simplicial maps, one can further improve the construction by using the notion of  annotation~\cite{BusaryevCabelloChen2012}, which are binary vectors  assigned to simplices.
This improved construction can be used to compute the persistence of the sparse zigzag Vietoris--Rips filtration of Sheehy~\cite{Sheehy2012}, without straightening out the zigzag. 
Whereas the approach by Sheehy~\cite{Sheehy2012} allows points to be deleted with a weighting scheme, Dey et al. provided an alternative method of sparsifying a Vietoris--Rips filtration via a more aggressive subsampling by collapsing input points in batches, thus  avoiding the weighting scheme completely. 
Their construction produces a sparsified filtration that is a $(3\log (1+\varepsilon)/2)$-approximation of the original filtration (see~\cref{sec:persistence}) and can be computed efficiently, where $0<\varepsilon \leq 1$ is a user-defined parameter.  

As a follow-up, Dey at al.~\cite{DeyShiWang2019} further provided a thorough computational comparison of the methods in \cite{Sheehy2012} and \cite{DeyFanWang2014}. They observed that, despite the linear size guarantee, the approach by Sheehy still ends up producing very large complexes due to the union step, and becomes unfeasible for high-dimensional data. On the other hand, the batch collapse from Dey et al.~\cite{DeyFanWang2014} is more space efficient, but still becomes prohibitively large for high-dimensional data. Inspired by their experimental findings, Dey at al.~\cite{DeyShiWang2019} proposed a new algorithm called SimBa, which obtains sparsified Vietoris--Rips filtrations through \emph{simplicial batch collapse}. 
Comparing with the approach in~\cite{DeyFanWang2014}, the key difference is that, instead of collapsing vertices based on their pairwise distance, SimBa builds another filtration by collapsing w.r.t. a \emph{set distance}. For two sets of points $A, B \subset P$, this set distance is defined as $d(A,B) = \min_{a \in B, \ b \in B}d(a,b)$. Although this may not seem a major modification, in practice, it ends up reporting significant improvements in memory and runtime. 

We have discussed finding cores only of Vietoris--Rips filtrations as sparsified filtrations. 
Kerber and Sharathkumar~\cite{KerberSharathkumar2013} extended the above ideas to sparsify \v{C}ech filtrations in two constructions.   
The first construction yields, for a fixed homological dimension, an approximate \v{C}ech filtration of linear size---just like in the Vietoris--Rips case---which is $(1+\varepsilon)$-interleaved with the original filtration. The key technical ingredient for this construction is a generalization to higher dimension of the \emph{well-separated pair decomposition} (WSPD), resulting in a \emph{well-separated simplicial decomposition} (WSSD). 
A WSSD decomposes a given point set $X$ into $O(n/\varepsilon^d)$ tuples; each $p$-tuple can be viewed as $p$ clusters of points satisfying that, if a ball contains at least one point of each cluster, a small expansion of the ball will contain all points~\cite{KerberSharathkumar2013}. 
In addition, for any $p$-simplex that can be formed with the points in $S$, there is a decomposition in $p+1$ clusters of $S$ where each cluster contains one vertex of the $p$-simplex. A drawback of this construction is that the constant in the size of the filtration depends exponentially on the dimension of the ambient space of the point cloud, whereas for the Vietoris--Rips construction, it depends on its doubling constant. Intuitively, the \emph{doubling constant} of a metric space is the minimum number of metric balls of some radius needed to cover any ball of double that radius. 

For the second construction, Kerber and Sharathkumar~\cite{KerberSharathkumar2013} generalized the Vietoris--Rips Lemma \cite{EdelsbrunnerHarer2022}, stating that the Vietoris--Rips complex at scale $r$ is contained in the \v{C}ech complex at scale $\sqrt{2}r$. 
This allows for the definition of a family of \emph{completion complexes}, verifying that the \v{C}ech complex at scale $r$ is contained in a completion complex at scale $(1+\varepsilon)r$. These complexes are parametrized by an integer $p$ and they are completely determined by their $p$-skeleton---of size at most $O(n^p)$---meaning that any higher-dimensional simplex is obtained combinatorially from the $p$-simplices. For a specific choice of $p$, one obtains a $(1+\varepsilon)$-interleaving with the original \v{C}ech filtration; see~\cite[Theorem 4]{KerberSharathkumar2013} for a precise statement. In this case, the construction does not depend on the dimension of the ambient space. For the proof of the approximation to the \v{C}ech complex, the key point is using coresets of the minimum enclosing ball of the point cloud. Here, the \emph{coreset} of a point cloud is a subset from which one can obtain an approximation of the minimum enclosing ball. 

M\'{e}moli and Okutan~\cite{MemoliOkutan2021} continued the above line of work and studied the simplification of \emph{filtered simplicial complexes}. 
A filtered simplicial complex indexed over $I \subset \Rspace$ is a family of simplicial complexes $(K_t)_{t \in I}$ such that $K_s$ is contained in $K_t$ ($K_s \hookrightarrow K_t$) for all $s \leq t$.   
It arises from the family of simplicial complexes that forms a filtration in persistent homology computation (see~\cref{sec:persistence}), such as the Vietoris--Rips or the \v{C}ech complexes.   
The authors aimed to simplify a filtered simplicial complex (i.e., reducing the number of simplices) while retaining its persistent homology. 
A key idea is the introduction of \emph{codensity}, defined on the vertices of a filtered simplicial complex, which quantifies the contribution of each vertex to the persistent homology (see~\cite{MemoliOkutan2021} for its technical definition based on vertex quasi-distance).
The authors showed that removing a vertex with zero codensity does not affect the homology. 
If one considers the Vietoris--Rips complex of a finite metric space, then the codensity of a point equals its distance to the nearest neighbour~\cite[Remark 4.1]{MemoliOkutan2021}. 
Similar to previous approaches~\cite{Sheehy2012,CavannaJahanseirSheehy2015}, the algorithm produces a simplified filtered complex by removing a vertex and all simplices containing it in an iterative way (i.e.,~based on codensity). 
The authors further introduced the notion of a \emph{simple filtered simplicial complex}, which is the one whose vertices all have positive codensity. 
They showed the existence of a unique (up to isomorphism) simple filtered simplicial complex (referred to as the \emph{core}) for any filtered simplicial complex such that their interleaving type distance is zero.  

Another way to obtain sparsified filtrations is via sequences of simplicial collapses. 
Bauer and Edelsbrunner~\cite{BauerEdelsbrunner2017} proved that, for a fixed parameter $r$, there is a collapsing sequence between the \v{C}ech, \v{C}ech-Delaunay, Delaunay (alpha), and the Wrap complexes. 
Since these complexes are homotopy equivalent, we may consider a collapsed complex as the core of the original complex by preserving both homological and homotopic properties.  
In particular, they proved the following theorem. 
\begin{theorem}[\v{C}ech-Delaunay Collapsing Theorem,~\cite{BauerEdelsbrunner2017}]
Let $X$ be a finite set of possibly weighted points in general position in $\Rspace^n$. Then, for every $r \in \Rspace$, 
\[
\Cech_r(X) \searrow \DelCech_r(X) \searrow \Del_r(X) \searrow \Wrap_r(X).
\]
\end{theorem}
This theorem implies that filtrations built from these four complexes have isomorphic persistent homology~\cite[Corollary 6.1]{BauerEdelsbrunner2017}. 

\subsection{Subsampling Methods Approximating Persistent Homology}
\label{sec:sc-cores-sampling}

We now discuss sparsification of simplicial complexes based on subsampling or bootstrapping. 
Cao and Monod~\cite{CaoMonod2022} introduced a bootstrapping technique to estimate the persistence diagram of a large input point cloud, whose direct persistence computation would be otherwise intractable. The authors proposed the following algorithm: given a finite subset $P$ of a metric space $(X, d)$ with a large number of points $|P|=N$, the algorithm takes a number of subsamples of $X$, each with $n \ll N$ points, and computes the persistence diagram of the Vietoris-Rips filtration of each subsample. 
The authors showed that the mean of the persistence diagrams of these bootstrap 
subsamples is a good approximation of the persistence diagram of the large data.  
Here, the \textit{mean} of persistence diagrams refers to the \textit{mean persistence measure}, obtained by first turning the persistence diagram into a discrete measure, where each point is substituted by a Dirac measure weighted by the multiplicity of the point in the diagram \cite{divol2019density,divol2021estimation}, and then taking the mean of the measures. 
Building on Divol and Lacombe's work on the convergence of the sample mean to the population mean of persistence diagrams \cite{divol2021estimation}, Cao and Monod derived explicit expressions for the approximation error of the convergence of the mean persistence measure through a bias-variance decomposition, justifying their bootstrap approach.

In a similar vein, G\'{o}mez and M\'{e}moli~\cite{GomezMemoli2024} studied Vietoris--Rips filtrations of subsamples. 
Given a compact metric space $(X,d)$, for a fixed dimension $p$ and an integer $n \geq 1$, they considered the persistence diagrams of $p$-dimensional Vietoris--Rips filtrations of all subsets of $(X,d)$ with cardinality at most $n \in \Nspace$. 
They worked with the notion of a \emph{curvature set} of a compact metric space $(X,d)$, first introduced by Gromov~\cite{Gromov1999}. 
For a fixed $n \in \Nspace$, the $n$-th curvature set of $X$ is a collection of \emph{all} $n \times n$ submatrices of the distance matrix of points in $X$ (with possible repetitions)~\cite{Gromov1999}. 
Starting from $(X,d)$, the algorithm takes subsamples of the distance matrix of $X$, and applies (Vietoris--Rips) persistent homology to each subsample, and aggregates the persistence diagrams of the subsamples by overlaying them into a single set of axes. The authors showed that the aggregated diagram is easy to compute, stable, and enjoys good discriminating power in classification tasks.    
Solomon et al.~\cite{SolomonWagnerBendich2022} also explored the idea of \emph{distributed persistence}, and demonstrated that the collection of persistence diagrams of \emph{many} small subsets of $(X,d)$ serves as a better invariant than a single persistence diagram of the entire space $(X,d)$.  

Although these approaches do not give rise to a single sparsified filtration, we may argue that the collection of subsamples and their corresponding filtrations serve as a replacement of the original data  filtration, and approximately preserve its persistence diagram.  
It remains open to develop concrete, geometric representations of such a collection.

\subsection{Sparsification Methods Approximating Spectral Properties}
\label{sec:sc-cores-spectral}

As an extension of spectral graph sparsification by Spielman and Srivastava~\cite{SpielmanSrivastava2011}, 
Osting et al.~\cite{OstingPalandeWang2020} introduced a subsampling algorithm for the $p$-skeleton of a simplicial complex that approximately preserves its $p$-dimensional up Laplacian. 
It relies on the effective resistance of simplices, generalizing \cref{eq:effective_resistances} using incidence matrix between simplices of adjacent dimensions. 
Given a simplicial complex $K$, and a fixed dimension $p \leq \dime(K)$, the algorithm runs as follows. 
It fixes all the simplices in $K$ up to dimension $p-1$ and then samples a number $q$ of $p$-simplices independently with replacement,  according to a probability distribution defined by their effective resistance. 
That way, the algorithm produces a sparsified simplicial complex $J$. 
Fixing a parameter $1/\sqrt{n_{p-1}}< \varepsilon \leq 1$, where $n_{p-1}$ is the number of $p-1$ simplices in $K$, for a sufficiently large value of $n_{p-1}$ and a precise value of $q$ (see \cite[Theorem 3.1]{OstingPalandeWang2020} for the full statement), Osting et al.~proved  that with probability at least $1/2$, they have, 
\[
(1-\varepsilon)  \mathbf{x}^T L_K \mathbf{x} \leq \mathbf{x}^T L_{J} \mathbf{x} \leq (1+\varepsilon) \mathbf{x}^T L_{K} \mathbf{x}, 
\]
for all $\mathbf{x} \in \Rspace^{n_{p-1}}$, where $L_K$ and $L_J$ are  the up Laplacian matrices for the original and sparsified simplicial complexes.

\subsection{Spectral Clustering of Simplicial Complexes Preserving Homology}
\label{sec:sc-cores-clustering}

Similar to the case of graphs (\cref{sec:graph-cores}), certain clustering methods of simplicial complexes may be considered as core finding as they preserve geometric or topological information of the complexes. 
Inspired by spectral clustering of graphs, Ebli and Spreemann~\cite{EbliSpreemann2019} introduced a spectral clustering of simplicial complexes. 
Given a $p$-dimensional simplicial complex $K$ of \emph{low homological complexity} (that is, $\beta_p(K) := \rank(\Hgroup_{p}(K)) \leq 10$), the algorithm produces a clustering of the simplices in $K$ of a fixed degree. 
The key idea is defining a \emph{harmonic embedding} $h$ of $K_p$, the $p$-skeleton of $K$, that is, 
\[
\varphi: K_p \to \Rspace^{\beta_p}, \varphi = \xi \circ j \circ i,  
\]
where $i: K_p \hookrightarrow \Cgroup_p(K)$ is the inclusion of $p$-skeleton $K_p$ to the $p$-chain group $\Cgroup_p(K)$, $j: \Cgroup_p(K) \to \Hcal_p(K)$ is the orthogonal projection onto the harmonic group $\Hcal_p(K)$ (i.e.,~the kernel of the Hodge Laplacian isomorphic to the homology group of $K$, $\Hgroup_p(K)$), and $\xi: \Hcal_p(K) \to \Rspace^{\beta_p(K)}$ is a choice of basis.  
In practice, Ebli and Spreemann simply chose an orthonormal basis  $h_1, \dots, h_{\beta_p(K)}$ for $\Hcal_p(K)$ and defined 
\[
\varphi(\sigma) = \left(\langle \sigma, h_1 \rangle_p, \dots, \langle \sigma, h_{\beta_p(K)} \rangle_p \right).
\]
In a lower-dimensional subspace of $\Rspace^{\beta_p(K)}$, the algorithm performs a subspace clustering method (e.g., independent component analysis \cite{HyvarinenOja2000}) of $\image(\varphi)$ and treats the output as the clustering of the original simplices.  
However, there are no theoretical guarantees for the above algorithm. The authors experimentally demonstrated that the simplices assigned to each cluster tend to reflect the presence of homology generators in the chosen dimension. For instance, when clustering the edges, they observed that edge clusters tend to respect the independent loops of the underlying space. The authors considered their algorithm to be complementary to the one presented by Osting et al.~\cite{OstingPalandeWang2020}.

\section{Future Research Opportunities}
\label{sec:future}

\subsection{Cores of Graphs}
\label{sec:future-graphs}

An aspect common to many methods using discrete curvatures to find graph cores through sampling is the lack of theoretical guarantees. Numerous experimental studies qualitatively show that these sampling techniques produce graphs similar to the originals and retain key structures. However, rigorous quantitative descriptions are scarce.

In practice, the primary sampling technique for sparsifying graphs using Forman-Ricci curvature involves simple thresholding to retain a certain percentage of edges. Several methods have been implemented to determine such a threshold, but it remains unclear which approach is optimal. A first direction for future research could be comparing all the methods proposed and evaluating in which scenarios they are more useful. In addition, future research could explore whether different sampling techniques, perhaps using the empirical distributions of edge curvatures, result in graph cores with distinct properties or impact the preservation of original graph structures.

As already mentioned, the Ollivier-Ricci curvature is mainly used for clustering. 
One reason it works for a graph with positive curvature can be explained by the inequalities in \cref{eq:bound-gap}, which ensure that it bounds the spectral gap of the Laplacian.
However, as reviewed in \cref{sec:graph-core-OR}, the edges representing the bottlenecks between clusters are those with the lowest negative curvature. A question one may ask is if and how this negative curvature also relates to the spectral gap.

Samal et al.~\cite{SamalSreejithGu2018} performed a comparative analysis of two discrete Ricci curvatures, the Forman-Ricci curvature and the Ollivier-Ricci curvature. 
They empirically showed a correlation between these two measures, which increases using the augmented Forman-Ricci curvature. However, the theoretical reasons for such a correlation require further study, which may shed some light on extracting graph cores with discrete curvatures. 

Finally, resistance curvature~\cite{DevriendtLambiotte2022,DevriendtOttoliniSteinerberger2024}, as a newer notion of discrete curvature for graphs, may be employed for finding graph cores, similar to the ways where Forman-Ricci the Ollivier-Ricci curvatures are employed.   

\subsection{Cores of Hypergraphs}
\label{sec:future-hypergraphs}

We have primarily reviewed methods for finding cores of hypergraphs based on sparsification and percolation techniques. However, given that hypergraphs inherently possess a complex underlying geometry, an outstanding gap in this area is the utilization of geometric methods, such as curvature-based approaches (e.g., \cite{LealRestrepoStadler2021}), for identifying the cores of hypergraphs.  

The notion of effective resistance was originally defined for graphs by Spielman and Srivastava~\cite{SpielmanSrivastava2011}, and then generalized to simplicial complexes by Osting et al.~\cite{OstingPalandeWang2020}, and more recently to hypergraphs by Aghdaei and Feng~\cite{AghdaeiFeng2022}.  
There are two possible future venues. First, whereas effective resistance has been used as a sampling criterion for sparsifying graphs~\cite{SpielmanSrivastava2011} and simplicial complexes~\cite{OstingPalandeWang2020}, its generalization to the sparsification of hypergraphs remains underexplored.  
Second, it would be interesting to define resistance curvature for hypergraphs and obtain curvature-based hypergraph cores. 
To define such a notion, we would need to define effective resistance of hypergraphs. Despite recent progress~\cite{AghdaeiFeng2022}, the question is: does there exist a canonical notion of effective resistance for hypergraphs, defined using an effective resistance matrix?   
In other words, is there an effective resistance matrix $R$ for hypergraph, in the form of $R := B L^{+} B^\top$, similar to \cref{eq:effective_resistances}, where $B$ is a boundary matrix for hypergraphs and $L$ is a hypergraph Laplacian?  
Whereas multiple generalizations from graph Laplacian to hypergraph Laplacian have appeared in the literature~\cite{AgarwalBransonBelongie2006,HeinSetzerJost2013,LiMilenkovic2018,SaitoMandicSuzuki2018,SaitoHerbster2023}, to the best of our knowledge, there does not appear to be a consensus as to the canonical notion of a hypergraph Laplacian.  
For instance, Rodr\'{\i}guez~\cite{Rodriguez2002,Rodriguez2003} introduced a version of the Laplacian matrix of a hypergraph based on the \emph{Laplacian degree} of a vertex and used it to obtain spectral-like results on partition problems in hypergraphs~\cite{Rodriguez2009}. 
Agarwal et al.~\cite{AgarwalBransonBelongie2006} defined higher order Laplacians as operators that measure variations on functions defined on $p$-chains formed by the vertex set. 
Hein et al.~\cite{HeinSetzerJost2013} made a connection between $p$-Laplacians of homogeneous hypergraphs with the total variation. 
Aktas and Akbas~\cite{AktasAkbas2022} introduced hypergraph Laplacians inspired by the simplicial Laplacian. 
It would be interesting to explore how various notions of resistance curvature (based on different Laplacians) affect the properties of hypergraph cores.  

Another potential direction for future research is the sparsification of hypergraphs while preserving their homological properties. Given the existence of various homology theories for hypergraphs (see~\cite{GasparovicPurvineSazdanovic2024} for a survey), we anticipate the development of distinct frameworks tailored to preserving different notions of homology. 

\subsection{Cores of Simplicial Complexes}
\label{sec:future-sc}

Within the area of finding the cores of simplicial complexes, there is still room for improvement in sparsification methods for filtrations approximately preserving the persistence diagram. Some ideas could be trying out different schemes for the batch collapse proposed by Dey et al. \cite{DeyFanWang2014} or trying to introduce some probabilistic techniques, useful in the sparsification techniques preserving spectral properties \cite{OstingPalandeWang2020}. In addition, there are other typical filtration constructions, known to have fewer simplices, that still pose computational constraints due to their size, such as filtrations based on the alpha complex. It could be interesting to see how the techniques described in  \cref{sec:sc-cores-PH} extend to this setting, which is already favorable for the purpose of reducing the size of the final complex.

On the other hand, concerning the method for sparsifying simplicial complexes preserving spectral properties described by Osting et al. \cite{OstingPalandeWang2020}, there are a few directions open for future research. First, from the condition obtained relating the Laplacians of the original and the sparsified complexes, it is not straightforward to see how the eigenvalues of the sparsified complex actually change. Finding a more direct relation in this context would be useful and illuminating. Additionally, the method allows sparsification only for simplices of dimension $p$, keeping the entire skeleton of the complex fixed up to dimension $p-1$.
 A potential direction of future research, already mentioned in \cite{OstingPalandeWang2020}, would be developing a methodology to sparsify across dimensions, without the need to fix the simplices in any dimension, and obtaining theoretical guarantees for it.

Finally, as already mentioned in \cref{sec:sc-cores-clustering}, a clear open question regarding the harmonic clustering algorithm described in \cite{EbliSpreemann2019} involves the derivation of theoretical guarantees for this method.
\section{Conclusion}
\label{sec:conclusion}

In this paper, we survey geometric and topological methods to extract the cores of graphs, hypergraphs, and simplicial complexes. 
After studying the collection of recent works, we found that there are noticeable imbalances among core-finding methods for higher graphs, revealing research gaps and hence opportunities. 
While many core-finding methods for graphs focus on sampling/thresholding with discrete curvatures (\cref{table:graph-core-summary}), similar methods are underdeveloped for hypergraphs and simplicial complexes with theoretical guarantees.  
Whereas there are a few recent developments in percolation theory for hypergraphs (\cref{table:hypergraph-core-summary}), there is plenty of room to grow. The percolation theory for simplicial complexes does not yet lead to the computation of well-defined cores.  
Sparse filtrations of simplicial complexes have been developed for preserving persistent homology (\cref{table:sc-core-summary}), so the question is, can one obtain sparsifications that preserve homological and spectral properties simultaneously? 
Spectral sparsification methods have been extended from graphs to hypergraphs and to simplicial complexes; however, the extracted cores may vary depending on different notions of Laplacians for higher graphs. 
We aim to use this survey (cf.~\cref{table:graph-core-summary,table:hypergraph-core-summary,table:sc-core-summary}) to inspire new research on the simplification of higher graphs.

\section*{Acknowledgments}
We are grateful to the 3rd Women in Computational Topology (WinCompTop) workshop for initiating our research collaboration. 
BW was supported in part by grants from the NSF (DMS-2301361 and IIS-2205418) and DOE (DE-SC0021015). 
SP was supported in part by grants from the NSF (CCF-2142713) and the National Cancer Institute of the National Institutes of Health under Award Number U54CA272167. The content is solely the responsibility of the authors and does not necessarily represent the official views of the National Institutes of Health. 
CL carried out this work under the auspices of INdAM-GNSAGA and partially within the activities of ARCES (University of Bologna). 
IGR is funded by a London School of Geometry and Number Theory--Imperial College London PhD studentship, which is supported by the Engineering and Physical Sciences Research Council [EP/S021590/1]. 
AS was supported in part by the European Research Council under Grant CoG 2015-682172NETS, within the Seventh European Union Framework Program; by the Swiss Government Excellence Scholarship; and by the Wallenberg AI, Autonomous Systems and Software Program (WASP), funded by the Knut and Alice Wallenberg Foundation.
LZ was supported in part by the Institute for Computational and Experimental Research in Mathematics at Brown University. 

We would like to thank Marzieh Eidi for discussions on the properties and challenges of the Ollivier-Ricci curvature. We would like to thank Guanqun Ma for generating the surfaces shown in \cref{fig:OR}. We would like to thank Heather A. Harrington and Polina Turishcheva for preliminary discussions on this project.

\bibliographystyle{plain}
\bibliography{refs-survey}

\end{document}